\documentclass[moor]{informs1}              

\usepackage{tikz}
\usetikzlibrary{shapes}
\newcommand*\circled[1]{\tikz[baseline=(char.base)]{
        \node[shape=rounded rectangle,draw,inner sep=2pt] (char) {#1};}}
\usepackage{amssymb}
\def\vomega{\mbox{\boldmath $\omega $}}
\usepackage{bbm}
\newcommand{\Mod}[1]{\mathrm{mod}\ #1}

\usepackage{natbib}
 \NatBibNumeric
 \def\bibfont{\small}%
 \def\BIBand{and}%
 \bibpunct[, ]{[}{]}{,}{n}{}{,}%

\usepackage[colorlinks=true,breaklinks=true,bookmarks=true,urlcolor=blue,
     citecolor=blue,linkcolor=blue,bookmarksopen=false,draft=false]{hyperref}

\def\EMAIL#1{\href{mailto:#1}{#1}}

\TheoremsNumberedThrough     

\EquationsNumberedThrough    

\begin{document}


 \RUNAUTHOR{Fibich and Levin}

 \RUNTITLE{Funnel Theorems for Spreading on Networks}

 \TITLE{Funnel Theorems for Spreading on Networks}

\ARTICLEAUTHORS{%
\AUTHOR{Gadi Fibich, Tomer Levin}
\AFF{Department of Applied Mathematics, Tel Aviv University, Tel Aviv 6997801, Israel, \EMAIL{fibich@tau.ac.il}, \EMAIL{levintmr@gmail.com}}
} 

\ABSTRACT{%
We derive novel analytic tools for the discrete Bass model, which models the diffusion of new products on networks.
We prove that the probability that any two nodes adopt by time~$t$, is greater than or equal to the product of the probabilities that each of the two nodes adopts by time~$t$. We introduce the notion of an ``influential node'', and use it  to determine whether the above inequality is strict or an equality. We then use the above inequality to prove the ``funnel inequality'', which relates the adoption probability of a node to the product of its adoption probability on two sub-networks. We introduce the notion of a ``funnel node'', and use it to determine whether the funnel inequality is strict or an equality.   The above analytic tools can be exptended to epidemiological models on networks. 
  
    We then use the funnel theorems to derive a new inequality for diffusion on circles and a new explicit expression for the adoption probabilities of nodes on two-sided line, and to prove that the adoption level on one-sided lines is strictly slower than on anisotropic two-sided lines, and that the adoption level on multi-dimensional Cartesian networks is bounded from below by that on one-dimensional networks.
}%


\MSCCLASS{Primary: 92D25; secondary: 91B99, 90B60}

\maketitle

%
\section{Introduction.}
%
Diffusion of new products is a classical problem in marketing~\cite{Mahajan-93}.
The diffusion starts when the product is first introduced into the market,
and progresses as more and more people adopt the product.
The first mathematical model of diffusion of new products
was introduced by Bass~\cite{Bass-69}.  In this model, individuals adopt a new product because
of {\em external influences} by mass media and {\em internal influences} (peer effect, word-of-moth)
by individuals who have already adopted the product.
This seminal study inspired a
huge body of theoretical and empirical research~\cite{Hopp-04}.

The Bass model, as well as most of this follow-up research, were carried out using compartmental models, which are typically given by deterministic  ordinary differential equations. Such models implicitly assume that
all individuals within the population are equally likely to influence each other, i.e., that the underlying social network is a homogeneous complete graph.
In more recent years, 
research on diffusion of new products gradually shifted 
to discrete Bass models on networks, in which the 
adoption decision of each individual is stochastic. 
The discrete Bass model allowed for heterogeneity among individuals,  
and for implementing a social network structure, whereby individuals are only influenced by adopters who are also their peers. 


Initially, discrete Bass models on networks were studied numerically, 
using agent-based simulations (see, e.g.~\cite{Garber-04,Goldenberg-09,GLM-02}).
To {\em analytically compute} the adoption probabilities of nodes in discrete Bass models on networks, one has to start from the    
 {\em master (Kolmogorov) equations} for the Bass model, which are~$2^M-1$ coupled linear ODEs, where~$M$ is the number of nodes (see, e.g., \cite[Section 3.1]{fibich2022diffusion}. Therefore, in order to explicitly solve these equations, one needs to reduce the number of ODEs significantly. 
 
    At present, there are two analytic techniques for solving the master equations explicitly, without making any ``mean-field'' type approximation. The first is based on {\em utilizing symmetries} of the master equations, in order to reduce the number of equations. This approach was  applied to homogeneous circles~\cite{OR-10}
    and to homogeneous and inhomogeneous complete networks~\cite{DCDS-22}.   
    The second approach is based on the {\em indifference principle}~\cite{Bass-boundary-18}. This analytic tool
    simplifies the explicit calculation of adoption probabilities, by replacing the original network with a simpler one.   
    The indifference principle has been used to compute the adoption probabilities of nodes on  bounded and unbounded lines, on circles, and on percolated lines~\cite{fibich2020percolation,Bass-boundary-18}.
   
\begin{figure}[ht!]
\begin{center}
\scalebox{1}{\includegraphics{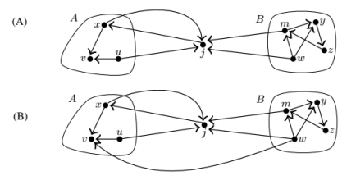}}
\caption{Illustration of funnel nodes. (A) Node $j$ is a vertex cut between $A:=\{x,v,u\}$ and~$B:=\{m,w,y,z\}$, and so is also a funnel node. (B)~Node~$j$ is not a vertex cut between~$A$ and~$B$. Nevertheless,  since there is no node in $A\cup B$ which is influential to $j$ both in~${\cal N}^A$ and in~${\cal N}^B$, node~$j$ is a funnel node of~$A$ and $B$. }
\label{fig:funnel}
\end{center}
\end{figure}

   In this paper, we introduce a third technique - the ``funnel theorems''
   (Section~\ref{sec:funnel-thms}). 
%
Choose some node~$j$, 
and divide the remaining nodes 
into two subsets of nodes: $A$ and $B$ (see Figure~\ref{fig:funnel}).
The funnel theorems provide the relation between the adoption probability of node~$j$ in the original network, with the product of the adoption probability of~$j$ on the two sub-networks~$\{j,A\}$ and~$\{j,B\}$,
which in many cases is easier to compute. The funnel relation is an equality if~$j$ is a {\em vertex cut}
(see Figure~\ref{fig:funnel}A), or more generally  if~$j$ is a~{\em funnel node}
(see Figure~\ref{fig:funnel}B), and is otherwise
 a strict inequality.

%

To prove the funnel theorems, we first prove that 
the probability any pair of nodes to be adopters is greater or equal 
than the product of the adoption probabilities of each of the nodes
(Section~\ref{sec:Si,j>-SiSj}). In other words, the correlation of adoption at the same time between any pair of nodes is non-negative.
This inequality is not only of interest by itself, but is also an important analytic tool. For example, we recently used it in~\cite{bounds-23} to derive optimal lower and upper bounds for the adoption level on any network.

To illustrate the power of the funnel theorems, we apply them to circular and Cartesian networks. Thus, 
\begin{enumerate}
	\item We derive a novel inequality for diffusion on circles (Theorem~\ref{thm:alpha_q_circle}).
	\item 
We derive a new explicit expression for the adoption probability of nodes on isotropic and nonisotropic two-sided lines (Theorem~\ref{thm:fj_twosided_new}).

\item We prove that the adoption level on the one-sided line is strictly slower than on two-sided isotropic and anisotropic lines
(Theorem~\ref{thm:f_line^one-sided<f_line_two-sided}). This improves on~\cite{Bass-boundary-18} by extending the proof to the anisotropic case.
In addition, the new proof is considerably simpler.

\item We prove that the adoption level on infinite multi-dimensional one-sided and two-sided Cartesian networks is strictly higher than on the infinite line    
(Theorem~\ref{thm:f_1D_bound}). 
\end{enumerate}

Finally, we note that our results are also relevant to {\em spreading of epidemics on networks}.  We discuss this in Section~\ref{sec:epidemiological}, and compare our results with those in the epidemiological literature.

\section{Theory review.}
\label{sec:main}
\subsection{Discrete Bass model.}
	\label{sec:discrete_Bass}
	We begin by introducing the discrete Bass model. A new product is introduced at time $t=0$ to a network with $M$ nodes, denoted by~${\cal M}:=\{1,\ldots,M\}$. Let~$X_j(t)$ denote the {\em state} of node/individual~$j$ at time~$t$, so that 
	\begin{equation*}
		X_j(t)=\begin{cases}
			1, \qquad {\rm if}\ j\ {\rm is\ an\ adopter\ of\ the\ product\ at\ time}\ t,\\
			0, \qquad {\rm otherwise.}
		\end{cases}
	\end{equation*} 
	Since all individuals are nonadopters at $t=0$,
	\begin{subequations}
	\label{eq:dbm}
	\begin{equation}
		\label{eq:general_initial}
		X_j(0)=0, \qquad j\in {\cal M}.
	\end{equation}
The adoption decision is irreversible, i.e., once a node adopts the product, it remains an adopter for all later times.

The adoption of nodes is stochastic, as follows. 
	Node~$j$ experiences {\em external influences} by mass media to adopt at the rate of~$p_j$. 
The underlying social network is represented by a  directed graph with positive weights, such that the weight of the edge from node~$k$ to node~$j$ 
is denoted by~$q_{k,j}> 0$, and $q_{k,j}=0$ if there is no edge from $k$ to $j$. Thus, if~$k$ already adopted the product and~$q_{k,j}>0$, its rate of {\em internal influence} on~$j$ to adopt is~$q_{k,j}$.
Since a nonadopter does not influence itself to adopt,
$$
q_{j,j}\equiv 0, \qquad j\in {\cal M}.
$$	
	Finally, internal and external influence rates are additive. Therefore,  
	{\em the adoption time~$\tau_j$ of~$j$ is a random variable, which is exponentially distributed} at the rate of
	\begin{equation}
		\label{eq:lambda_j}
		\lambda_j(t):=p_j+\sum_{k\in {\cal M}}q_{k,j}X_k(t), \qquad j\in {\cal M}, \quad t>0.
	\end{equation}
Thus, time is continuous, and the adoption rate of~$j$ increases whenever one of its peers becomes an adopter. 
\end{subequations}	

The maximal rate of internal influences that can be exerted on~$j$ (which is when all its neighbors/peers  are adopters) is denoted by
\begin{equation}
	\label{eq:q_on_node}
	q_j:=\sum_{\substack{k \in {\cal M}}}q_{k,j}.
\end{equation}
%
%
	The underlying network of the discrete Bass model~\eqref{eq:dbm} is denoted by
	$$
	{\cal N} =  {\cal N}({\cal M}, \{p_k\}_{k\in {\cal M}}, \{q_{k,j}\}_{k,j\in {\cal M}}).
	$$ 

Our  goal is to explicitly compute the adoption probability of nodes
\begin{equation*}
	f_j(t) := \mathbb{P}(X_j(t)=1)=\mathbb{E}\left[X_j(t)\right], \qquad j \in {\cal M},
\end{equation*}
and use this to compute the expected {\em fraction of adopters}
({\em  adoption level}) 
	\begin{equation*}
		\label{eq:number_to_fraction}
		f(t):=\frac{1}{M}\mathbb{E} \left[\sum_{j =1}^M X_j(t\right] 
		=\sum_{j =1}^Mf_j(t),
	\end{equation*}
where~$\sum_{j=1}^M X_j(t)$ is the number of adopters at time~$t$.
In most cases, it is easier to compute the corresponding nonadoption probabilities
\begin{equation*}
\label{eq:nonadoption_prob}
[S](t):=1-f(t), \qquad [S_j](t):=1-f_j(t)= \mathbb{P}(X_j(t)=0), \quad j \in {\cal M}.
\end{equation*}

	\subsection{Dominance and indifference principles.}
	\label{sec:dominance}
	
	The {\em dominance principle} is useful for comparing the adoption probabilities of nodes 
	in two networks.
	Let us begin with 
	\begin{definition}[network dominance]
		\label{def:dominance}
		Consider the discrete Bass model~\eqref{eq:dbm} on networks
$$
{\cal N}^A := N({\cal M}, \{p_k^A\}_{k=1}^M, \{q_{k,j}^{A}\}_{k,j=1}^M), \qquad  {\cal N}^B := N({\cal M}, \{p_k^B\}_{k=1}^M, \{q_{k,j}^{B}\}_{k,j=1}^M),
$$		
both with~$M$ nodes.
		We say that ``${\cal N}^A$ is dominated by ${\cal N}^B$'' and denote  ${\cal N}^A \preceq {\cal N}^B$,  if 
		$$
		p_j^{A} \leq p_j^B \quad \mbox{for all~j} \qquad \mbox{and} \qquad  q_{k,j}^{A} \leq q_{k,j}^{B} \quad \mbox{for all~} k \not=j.
		$$
		We say that  ``${\cal N}^A$ is strongly dominated by ${\cal N}^B$'' and denote   ${\cal N}^A \prec {\cal N}^B$, if at least one of these $M^2$ inequalities is strict.
	\end{definition}
	
	\begin{theorem}[dominance principle for nodes~\cite{Bass-boundary-18}]
	\label{lem:dominance-principle-f_j}
	Consider the discrete Bass model~\eqref{eq:dbm} on networks ${\cal N}^A$ and~${\cal N}^B$, both  with $M$~nodes.
	 If ${\cal N}^A \preceq {\cal N}^B$, then the adoption probability of any 
	 node in network~${\cal N}^A$ is lower than or equal to its adoption probability in network~${\cal N}^B$, i.e., 
	$$
	f_j^A(t) \leq f_j^B(t), \qquad 0 \le t < \infty, \qquad  j \in {\cal M}.
	$$
%
\end{theorem}

%

Let $\Omega\subset {\cal M}$ denote a subset of the nodes, and let
$$
[S_\Omega](t):= \mathbb{P}\left(X_j(t) = 0,   j\in \Omega \right)
$$
denote the probability that all nodes in~$\Omega$ did not adopt by time~$t$.
The {\em indifference principle} simplifies the explicit calculation of~$[S_\Omega]$, 
by replacing the original network with a  network with a modified edge  structure, such that the value of~$[S_\Omega]$ remains unchanged, but its explicit calculation is simpler. To do that, we need to 
be able to distinguish between edges that influence the nonadoption probability~$[S_\Omega](t)$ and those that do not.

\begin{definition}[influential and non-influential edges to~$\Omega$~\cite{Bass-boundary-18}]
	\label{def:influential-edge}
	Consider a directed network with  $M$~nodes \emph{(}if the network is undirected, replace each undirected edge by two directed edges\emph{)}.
	Let $\Omega \subsetneqq {\cal M}$ be a subset of the nodes, and let $\Omega^{\rm \emph{c}} := {\cal M} \setminus \Omega $ be its complement.

A directed edge \circled{$k$}~$\to$~\circled{$m$}  is called ``influential to $\Omega$", if the following two conditions hold:
\begin{enumerate}
	\item $k\in \Omega^{\emph{\text{c}}}$, and
	\item either $m \in \Omega$, or there is a finite sequence of directed edges from node~$m$ to some node $u \in \Omega$, which does not go through node~$k$.
	\end{enumerate}  

	A directed edge \circled{$k$}~$\to$~\circled{$m$} is called ``non-influential to $\Omega$'' if one of the following three conditions holds:
	\begin{enumerate}
		\item  $k \in \Omega$, or
		\item  $k  \in \Omega^{\rm \emph{c}}, m  \in \Omega^{\rm \emph{c}}$, and there is no finite sequence of directed edges from node~$m$ to $\Omega$, or
		\item  $k  \in \Omega^{\rm \emph{c}}, m  \in \Omega^{\rm \emph{c}}$, and all finite sequences of directed edges from node~$m$ 
		to~$\Omega$ go through the node~$k$.
	\end{enumerate}

Thus, any edge is either ``influential to~$\Omega$'' or  
``non-influential to $\Omega$".
\end{definition}

\begin{theorem}[Indifference principle~\cite{Bass-boundary-18}]
	\label{thm:Indifference}
	Consider the discrete Bass model \eqref{eq:dbm}, and 
	let $\Omega \subsetneqq {\cal M}$ be a subset of the nodes.
	 Then $[S_{\rm \Omega}](t)$
	remains unchanged if we remove or add edges which are non-influential to~$\Omega$.
\end{theorem}

\subsection{Strong dominance principle for nodes.}
\label{sec:sharp-dominance}

In Theorem~\ref{lem:dominance-principle-f_j} we saw that if 
${\cal N}^A \preceq {\cal N}^B$, then $f_j^A(t) \leq f_j^B(t)$  for any $j \in {\cal M}$.
In~\cite{Bass-boundary-18} it was also showed that
if ${\cal N}^A \prec {\cal N}^B$, then 
the adoption level in~${\cal N}^A$ is strictly lower than in~${\cal N}^B$, i.e.,  
$f^{A}(t) < f^{B}(t)$. Since $f^{A} = \frac1M \sum_{j \in \cal M} f_j^A$
and  $f^{B} = \frac1M \sum_{j \in \cal M} f_j^B$,
the condition ${\cal N}^A \prec {\cal N}^B$ implies that 
$f_j^A(t) < f_j^B(t)$ for at least one node. 
%
In order to fully characterize the nodes for which the condition ${\cal N}^A \prec {\cal N}^B$ implies that $f_j^A(t) < f_j^B(t)$, we introduce

 \begin {definition} [Influential node] \label{def:influential_nodes} 
Let $\Omega \subset \cal M$. 
We say that ``node~$m$ is influential to~$\Omega$'' if $m \in \Omega$, or if $m \in \Omega^c$ and there is finite sequence of directed edges from~$m$ to~$\Omega$.
\end{definition}

Thus,  $m \in \Omega^c$ in an influential node to~$\Omega$, if and only if 
there is an edge emanating from~$m$ which is influential to~$\Omega$.


\begin{lemma}[strong dominance principle for nodes]
	\label{lem:strong-dominance-principle-nodes}
	Consider the discrete  Bass model~\eqref{eq:dbm} on two networks ${\cal N}^A$ and~${\cal N}^B$, both with $M$~nodes.
If ${\cal N}^A \prec {\cal N}^B$, then  the adoption probability of any 
node in network~${\cal N}^A$ is strictly lower than its adoption probability in network~${\cal N}^B$, i.e., 
$$
f_j^A(t) <f_j^B(t), \qquad 0 \le t < \infty, 
$$
  if and only if  at least one of the following two conditions holds:
   \begin{enumerate}
	\item  $p_m^{A}< p_m^{B}$ and node~$m$ is influential to~$\Omega = \{j\}$.
	\item $q_{m,k}^{A}< q_{m,k}^{B}$ and the edge $\circled{m} \rightarrow \circled{k}$ is influential to~$\Omega = \{j\}$.
\end{enumerate}
\end{lemma}
\proof{Proof.}
	The proof is similar to that of Theorem~\ref{lem:dominance-principle-f_j} in \cite{Bass-boundary-18}. \Halmos
\endproof

\subsection{Multivariate Chebyshev's integral inequality.}

Let us recall the classical one-dimensional Chebyshev's integral inequality with weights:
\begin{lemma}
	\label{lem:chebyshev-1D-weights}
	Let~$f(x)$ and~$g(x)$ be both non-decreasing (or both non-increasing) functions in~$I:=\left[a,b\right]$, where infinite boundaries are allowed. Let~$w(x)$ be a positive function in~$I$, for which~$\int_a^b w(x)\,dx = 1$.
	 Then
	\begin{equation*}
		\int_a^b f\left(x\right)g\left(x\right)w(x)\,dx  \geq \int_0^1f\left(x\right)w(x)\,dx \int_0^1g\left(x\right)w(x)\,dx.
	\end{equation*} 
	Furthermore,  an equality holds  if and only if  either~$f$ or~$g$ are independent of~$x$.
\end{lemma}
\proof{Proof.}
	See Appendix~\ref{app:chebyshev-1D-weights}. \Halmos
\endproof
	
%
%


The multi-dimensional extension of Lemma~\ref{lem:chebyshev-1D-weights} with~$[a,b]=[0,1]$ and~$w(x)\equiv 1$ is 
\begin{lemma}[Chebyshev's multi-dimensional integral inequality]
	\label{lem:chebyshev-multi-D}
	Let~$f\left(x_1,\ldots,x_D\right)$ and $g\left(x_1,\ldots,x_D\right)$ be both non-decreasing (or both non-increasing) functions in~$[0,1]^D$ 
	with respect to each~$x_i$, where $i=1,\ldots, D$.
	Then 
	\begin{equation*}
		\begin{aligned}
			\int_{\left[0,1\right]^D}f&\left(x_1,\ldots,x_D\right)g\left(x_1,\ldots,x_D\right)dx_1\cdots dx_D \geq\\
			&\int_{\left[0,1\right]^D}f\left(x_1,\ldots,x_D\right)dx_1\cdots dx_D\int_{\left[0,1\right]^D}g\left(x_1,\ldots,x_D\right)dx_1\cdots dx_D.
		\end{aligned}
	\end{equation*}
	Furthermore,  an equality holds if and only if  for any $i=1, \dots, D$, 
	either~$f$ or~$g$ are independent of~$x_i$.
\end{lemma}
\proof{Proof.}
	See Appendix~\ref{app:chebyshev-multi-D}. \Halmos
\endproof

\section{ $[S_{i,j}] \ge [S_i]\, [S_j] $.}
\label{sec:Si,j>-SiSj}

Consider the adoption probability 
$$
f_{i,j} (t):=\mathbb{P}(X_i(t) = X_j(t) = 1)
$$
that both~$i$ and~$j$ are adopters.
If~$i$ adopts, it may influence node~$j$ to adopt as well. 
Therefore, if we know that~$i$ is an adopter, this increases the likelihood 
that~$j$ is also an adopter, i.e.,  
$
 (f_j  \,|\, f_i) \ge f_j.
$
Hence,  it is reasonable to expect that
$
f_{i,j} = f_i \, (f_j  \,| \, f_i)  \ge  f_i \, f_j.
$
To prove this inequality, it is convenient  to reformulate it using the nonadoption probability 
$$
[S_{i,j}](t):=\mathbb{P}(X_i(t) = X_j(t) = 0)
$$
that both~$i$ and~$j$ are nonadopters.

\begin{lemma}
	\label{lem:fi,j-fifj=Si,j-SiSj}
	$$
	[S_{i,j}]  -  [S_i]\, [S_j] = f_{i,j}  - f_i f_j.
	$$
\end{lemma}
\proof{Proof.}
	See Appendix~\ref{app:fi,j-fifj=Si,j-SiSj}. \Halmos
\endproof
Therefore, 
$$
f_{i,j}  \ge  f_i \, f_j \iff [S_{i,j}]  \ge  [S_i]\, [S_j].
$$
Indeed, we have the following result:
\begin{theorem}
	\label{thm:[Si,Sj>=[Si][Sj]}
	Consider the discrete Bass model~\eqref{eq:dbm}.
	Then for any two nodes $i,j \in {\cal M}$, 
	\begin{equation}
	\label{eq:[Si,Sj>=[Si][Sj]}
		[S_{i,j}](t) \ge [S_{i}](t) \, [S_j](t), \qquad t \ge 0.
	\end{equation}
\end{theorem}
\proof{Proof.}
	By~\eqref{eq:dbm}, the stochastic adoption of~$j\in {\cal M}$ in the time interval~$(t,t+\Delta t)$ as~$\Delta t\rightarrow 0$ is given by the conditional probability
	\begin{equation}
		\label{eq:general_model}
		\mathbb{P}(X_j(t+\Delta t)=1\,| \, {\bf X}(t))=\begin{cases}
			\hfill 1,\hfill &\qquad {\rm if}\ X_j(t)=1,\\
			\left(p_j+\sum\limits_{\substack{k \in {\cal M}}}q_{k,j}X_{k}(t)\right)\Delta t, &\qquad {\rm if}\ X_j(t)=0,
		\end{cases}
	\end{equation}
	where ${\bf X}(t) := \left(X_1(t),\ldots, X_M(t)\right)$ is the {\em state of the network} at time~$t$.
		 Let $\Delta t >0$,  $t_{n} := n\Delta t$, and $\vomega^{n} := (\omega_{1}^{n},\ldots,\omega_{M}^{n}) \in [0,1]^{M}$. 
We define the time-discrete realization
\begin{equation}
	\label{eq:tilde-Xj^A-def}
	\widetilde{X}_{j}(t_{n})
	:=X_{j}(t_{n};\{\vomega^{n}\}_{n=1}^\infty,\Delta t), \qquad  j \in {\cal M}, \quad n = 0,1, \dots, 
\end{equation}
 of~\eqref{eq:general_model} as follows:

\renewcommand\labelitemii{$\bullet$} 		 
\renewcommand\labelitemiii{$\bullet$} 
\renewcommand\labelitemiv{$\bullet$} 
\begin{itemize}
	\item $\widetilde{X}_{j}(0) := 0$ for $ j \in {\cal M} $
	\item for $n=1,2,\ldots$ \begin{itemize}
		\item for $ j \in {\cal M} $ \begin{itemize}
			\item if $\widetilde{X}_{j}(t_{n-1}) = 1$, then $\widetilde{X}_{j}(t_{n}) := 1$
			\item if $\widetilde{X}_{j}(t_{n-1}) = 0$, then
			\begin{itemize}
				\item if $\omega_{j}^{n} \leq \Big( p_j(t_{n}) + \sum_{k\not= j}q_{k,j}(t_{n}) \widetilde{X}_{k}(t_{n-1}) \Big)\Delta t$, then $\widetilde{X}_{j}(t_{n}) := 1$ 
				\item else $\widetilde{X}_{j}(t_{n}) := 0$ 
			\end{itemize}
		\end{itemize}
		\item end
	\end{itemize}
	\item end
\end{itemize} 
Let us fix $t>0$. Let $N \gg 1$ and~$\Delta t  := \frac{t}{N}$. 
Note that as $N \to \infty$, $\Delta t \to 0$ and~$t_N \equiv t$.
Therefore, to prove~\eqref{eq:[Si,Sj>=[Si][Sj]}
it is sufficient to show that for any $0 < \Delta t \ll 1$,
\begin{equation}
	\label{eq:[S_iS_j]-discrete}
	[S_{i,j}](t_N;\Delta t) \ge [S_i](t_N;\Delta t) [S_j](t_N;\Delta t).
\end{equation}
To do that, we first note that
\begin{align*}
&[S_i](t_N;\Delta t)
= \int_{[0,1]^{M \times N}} G_i
\, d\vomega^{1}\cdots d\vomega^{N},
\qquad
[S_j](t_N;\Delta t)
= \int_{[0,1]^{M \times N}}  G_j  
\, d\vomega^{1}\cdots d\vomega^{N},
\\
&[S_{i,j}](t_N;\Delta t) = \int_{[0,1]^{M \times N}}  G_i G_j \, d\vomega^{1}\cdots d\vomega^{N},
\end{align*}
where
$$
G_i = G_i(t_N;\{\vomega^{n}\}_{n=1}^N,\Delta t):= \mathbbm{1}_{X_{i}(t_{N};\{\vomega^{n}\}_{n=1}^N,\Delta t)=0},
$$
and similarly for~$G_j$.

We claim that the function~$G_i$ and~$G_j$ are both
non-decreasing in~$\left[0,1\right]^{M \times N}$ 
with respect to any~$\omega_m^n$, where
$m \in {\cal M}$ and $n=1, \dots, N$.
Therefore, inequality~\eqref{eq:[S_iS_j]-discrete}~follows from  
Chebyshev's multi-dimensional integral inequality (Lemma~\ref{lem:chebyshev-multi-D}).

To prove this claim, note that if we increase $\omega_{m_0}^{n_0}$ by a factor of~$\beta>1$, this is completely equivalent to decreasing $p_{m_0}(t_{n_0})$ and $\{q_{k,{m_0}}(t_{n_0})\}_{k \not=m_0}$ by $\beta$, since in the definition of $\widetilde{X}_j$, $\omega_{m}^{n}$ only appears in the condition
$\omega_{m}^{n} \leq \Big( p_m(t_{n}) + \sum_{k\not= m}q_{k,m}(t_{n}) \widetilde{X}_{k}(t_{n-1}) \Big)\Delta t$.
Therefore, from the proof of the dominance principle, see~\cite[eq.~(3.4)]{Bass-boundary-18}, for any~$n$, we have that  
either $\widetilde{X}_i\left(t_{n}\right)$ and~$\widetilde{X}_j\left(t_{n}\right)$
decrease, or they remain unchanged.
Hence, either~$G_i$ and~$G_j$  increase or they remain unchanged.  With this, the proof of inequality~\eqref{eq:[Si,Sj>=[Si][Sj]} concludes. \Halmos
\endproof

\subsection{When does $[S_{i,j}] > [S_i]\, [S_j]$?}

The condition for inequality~\eqref{eq:[Si,Sj>=[Si][Sj]} to be strict 
makes use of the notion of an {\em influential node} (Definition~\ref{def:influential_nodes}).


\begin{theorem}
\label{thm:[Si,Sj>[Si][Sj]}
Consider the discrete Bass model~\eqref{eq:dbm}.
\begin{enumerate}
		\item If there exists a node in~${\cal M}$ which is influential  to~$i$ and to~$j$, then 
	\begin{equation}
	\label{eq:[Si,Sj>[Si][Sj]}
	[S_{i,j}](t) > [S_{i}](t) \, [S_j](t), \qquad t > 0.
	\end{equation}
	\item If, however, there is no node which is influential to~$i$ and to~$j$, then 
\begin{equation}
	\label{eq:[Si,Sj=[Si][Sj]}
	[S_{i,j}](t) = [S_{i}](t) \, [S_j](t), \qquad t \ge 0.
\end{equation}	
\end{enumerate}
\end{theorem}
\proof{Proof.}
By Chebyshev's integral inequality (Lemma~\ref{lem:chebyshev-multi-D}),  inequality~\eqref{eq:[S_iS_j]-discrete}  is an equality if and only if~$G_i$ and~$G_j$ depend on different coordinates.
The function $G_i$ depends on~$\omega_k^n$ if and only if
node~$k$ is influential to node~$i$. Therefore, inequality~\eqref{eq:[S_iS_j]-discrete} is
an equality if and only if there is no node which is influential to both~$i$ and~$j$.
Therefore, we prove~\eqref{eq:[Si,Sj=[Si][Sj]}.
            To prove~\eqref{eq:[Si,Sj>[Si][Sj]}, however, we also need to show that inequality~\eqref{eq:[S_iS_j]-discrete} remains strict as $\Delta t \to 0$.
            
 To see that, assume that node $m \in {\cal M}$ is influential to~$i$ and to~$j$.  Denote by~$\tau_m$ the random variable given by the adoption time of~$m$.
Let $H(x) = \mathbb{P}(\tau_m \le x)$ denote the CDF of~$\tau_m$,  and let
$h = H'$ denote its density.
Then by the law of iterated expectations
\begin{equation}
\label{eq:[S_{i,j}]-tau_m}
[S_{i,j}](t) = \int_{\tau_m=0}^\infty \mathbb{P}(X_i(t) = 0, X_j(t) = 0\, | \, \tau_m) \, h(\tau_m)  d\tau_m.
\end{equation}
$$
[S_i](t) =
\int_{\tau_m=0}^\infty \mathbb{P}(X_i(t) = 0 \, | \, \tau_m)  \, h(\tau_m) d\tau_m,
\quad
[S_j](t) =
\int_{\tau_m=0}^\infty \mathbb{P}(X_j(t) = 0 \, | \, \tau_m)  \, h(\tau_m) \, d\tau_m.
$$
For a given $\tau_m \ge 0$, the adoption of any node in ${\cal M} \setminus \{m\}$ in the original network
is identical to its adoption in network~$\widetilde{\cal N}$ with nodes
$\widetilde{\cal M}:={\cal M} \setminus \{m\}$ and with weights
$\tilde{p}_j := p_j+ q_{m,j} \mathbbm{1}_{t \ge \tau_m}$ and
$\tilde{q}_{k,j} := q_{k,j}$ for $k,j \in \widetilde{\cal M}$.
Therefore, by inequality~\eqref{eq:[Si,Sj>=[Si][Sj]}, 
\begin{equation}
            \label{eq:claim-conditional_tau_m}
           \mathbb{P}(X_i(t) = 0, X_j(t) = 0\, | \, \tau_m) \ge
            \mathbb{P}(X_i(t) = 0 \, | \, \tau_m) \, \mathbb{P}(X_j(t) = 0 \, | \, \tau_m).
\end{equation}
Combining~\eqref{eq:[S_{i,j}]-tau_m} and~\eqref{eq:claim-conditional_tau_m}, we have
\begin{align*}
[S_{i,j}](t)
\ge \int_{\tau_m=0}^\infty \mathbb{P}(X_i(t) = 0 \, | \, \tau_m) \, \mathbb{P}(X_j(t) = 0 \, | \, \tau_m)  \, h(\tau_m)  d\tau_m.
\end{align*}

Since node $m \in {\cal M}$ is influential to~$i$ and to~$j$ and $\tilde{p}_j = p_j+ q_{m,j} \mathbbm{1}_{t \ge \tau_m}$ is monotonically-decreasing in~$\tau_m$, then
by the dominance principle, the two conditional probabilities on the right-hand-side are strictly monotonically-increasing in $\tau_m$ for $0 \le \tau_m \le t$ and are  monotonically-increasing
for $0 \le \tau_m < \infty$.
Therefore, by Chebyshev's integral inequality with weights
(Lemma~\ref{lem:chebyshev-1D-weights}),
\begin{align*}
            [S_{i,j}](t) >
            \int_{\tau_m=0}^\infty \mathbb{P}(X_i(t) = 0 \, | \, \tau_m)  \, h(\tau_m)  d\tau_m
\int_{\tau_m=0}^\infty  \mathbb{P}(X_j(t) = 0 \, | \, \tau_m)  \, h(\tau_m)  d\tau_m
 =     [S_i](t) [S_j](t),
\end{align*}
as needed. \Halmos
\endproof

Since any node is influential to itself, we have
\begin{corollary}
	Consider the discrete Bass model~\eqref{eq:dbm}.
	If node~$i$ is influential to node~$j$, then 
	\begin{equation*}
			[S_{i,j}](t) > [S_{i}](t) \, [S_j](t), \qquad t > 0.
\end{equation*}
\end{corollary}

\section{Funnel theorems.}
  \label{sec:funnel-thms}

The funnel theorems make use of the following partition of nodes:

\begin{definition}[partition of nodes]
	Let $j \in {\cal M}$, and $A,B \subset {\cal M}$. We say that $\{A,B,\{j\}\}$ is a partition of 
	$ \cal{M}$, if ${\cal M} = A\cup B \cup \{{j} \}$, and  if $A$, $B$, and $\{j\}$ are	nonempty and mutually disjoint. 
\end{definition}

The adoption of node~$j$ in network~${\cal N}$
may be the result of three distinct direct influences:
\begin{enumerate}
	\item Internal influences on~$j$ by edges that arrive from~$A$.
	\item Internal influences on~$j$ by edges that arrive from~$B$.
	\item External influences on~$j$.
\end{enumerate}

 In order to identify the specific influence that led  to the adoption of~$j$, we define three networks, on which $j$ can only adopt due
to one of these three influences:

\begin{definition} [Networks~${\cal N}^{A,B,p_j}$, ${\cal N}^A$, ${\cal N}^B$, ${\cal N}^{p_j}$]
	\label{def:networks}
	Consider the discrete Bass model~\eqref{eq:dbm} with network  
	$
	{\cal N} := {\cal N}({\cal M}, \{p_k\}_{k=1}^M, \{q_{k,j}\}_{k,j=1}^M)
	$.
     Let  $\{A,B,\{j\}\}$ be a partition of 
	$ \cal{M}$. 
	We define four different networks with respect to this partition: 
	\begin{enumerate}
		\item ${\cal N}^{A,B,p_j}:= {\cal N}$ is the original network.
		\item Network ${\cal N}^A$ is obtained from~${\cal N}$ by removing all influences on node~$j$, except for directed edges from~$A$ to~$j$. Thus, we cancel the external influences on node~$j$  by setting~$p_j=0$, and we remove all direct links from~$B$ to~$j$ by setting~$q_{m,j} = 0$ for all~$m\in B$.
		\item Network ${\cal N}^B$ is defined similarly.
		\item Network ${\cal N}^{p_j}$ is obtained from~${\cal N}$  by removing all internal influences on node~$j$, but retaining  the external influence~$p_j$. Thus, we remove all direct links to~$j$ by setting~$q_{m,j} = 0$ for all~$m\in A \cup B$. 
	\end{enumerate}
\end{definition}
We denote the state of node~$j$ in each of these four networks
by
$$
X^{U}_j(t):=X^{N^U}_j(t), \qquad U \in \{\{A,B,p_j\}, A, B, p_j\}. 
$$

The funnel theorems below compare the nonadoption probability of node~$j$ in the original network, with the product of the three nonadoption probabilities of~$j$ due to each of the three distinct influences. 

\begin{theorem}[Funnel inequality]
	\label{thm:funnel_node_inequality}
	Consider the discrete Bass model~\eqref{eq:dbm}.
	Let $\{A,B,\{j\}\}$ be a partition of $ \cal{M}$.
	Then
	\begin{equation}
		\label{eq:funnel_inequality}
		\mathbb{P} \left(X^{A,B,p_j}_j(t)=0\right)  \ge
		\mathbb{P} \left(X^{A}_j(t)=0\right)
		\mathbb{P} \left(X^{B}_j(t)=0\right)
		\mathbb{P} \left(X^{p_j}_j(t)=0\right),  \qquad t \ge 0,
	\end{equation}
	where 
	\begin{equation}
	\label{eq:funnel_p_j}
	\mathbb{P} \left(X^{p_j}_j(t)=0\right) = e^{-p_j t}.
	\end{equation}
\end{theorem}
\proof{Proof.}
	By the indifference principle, all edges that emanate from~$j$ are non-influential to~$j$. Since this holds for all the four probabilities in~\eqref{eq:funnel_inequality}, in what follows, we can assume that no edges emanate from~$j$. 
	
	In principle, we need to compute the four probabilities in~\eqref{eq:funnel_inequality} using the four different networks from Definition~\ref{def:networks}. 
	We can simplify the analysis, however, by considering only two networks (which are also ``quite similar''), as follows.  
		Given network~${\cal N}$, we define network~${\cal N}^+$ by ``splitting'' node $j$ into three nodes  ${j}_A$,  ${j}_B$, 
	and~${j}_p$, such that:
	\begin{enumerate}
		\item Node  ${j}_A$ inherits from $j$ all the (one-sided) edges from~$A$ to~$j$, i.e.,
		$$
		{q}_{k,{j}_A}^+=
		\begin{cases}
			{q}_{k,j}, & k \in A,\\
			0, & k  \in B,
		\end{cases}
		$$
		and has ${p}_{{j}_A}^+=0$.
		\item Node  ${j}_B$ is defined similarly.
		\item  Node ${j}_p$ inherits from $j$ the external influences ${p}_{{j}_p}^+=p_j$, but has no incoming edges from~$A$ and $B$, i.e., ${q}_{k,{j}_p}^+ \equiv 0$ for all $k \in A \cup B$. 
			
		\item Since  no edges emanate from~$j$ in network~${\cal N}$, no edges emanate from nodes ${j}_A$, ${j}_B$, and~${j}_p$  in network~${\cal N}^+$.
		\item The weights of all nodes but~$j$ and of the edges between these nodes are the same in both networks.
		
	\end{enumerate}
	Let ${X}^{+}_{k}(t)$ denote the state of node~$k$ in network ${\cal N}^+$.
	By construction, 
	\begin{subequations}
		\label{eq:Prob(X^j,A)=Prob(X^j,tilde-A)}
		\begin{equation}
			\mathbb{P} \left(X^{A}_j(t)=0\right) = 
			\mathbb{P} \left({X}^{+}_{{j}_A}(t)=0 \right) ,
			\qquad 
			\mathbb{P} \left(X^{B}_j(t)=0\right) = 
			\mathbb{P} \left({X}^{+}_{{j}_B}(t)=0 \right) ,
			\end{equation}
		and
		\begin{equation}
			\mathbb{P} \left(X^{p_j}_j(t)=0\right) =  \left({X}^{+}_{{j}_p}(t)=0\right) .
		\end{equation}
	\end{subequations}
	
	In Appendix~\ref{app:ProofX^jA} we will prove that 
	\begin{equation}
		\label{eq:Prob(X^j,A)}
		\mathbb{P} \left(X_j(t)=0\right) =
		\mathbb{P} \left({X}^{+}_{{j}_p}(t)
		={X}^{+}_{{j}_A}(t)
		={X}^{+}_{{j}_B}(t)
		=0\right),
	\end{equation}
	 where ${X}_{k}(t)$ denotes the state of node~$k$ in network ${\cal N}$.
	Since $ {j}_p$ is an isolated node in~${\cal N}^+$, its adoption is independent of that of ${j}_A$ and ${j}_B$, and so
	\begin{equation}
		\label{eq:Prob(X^j,A)-indifference}
		\mathbb{P} \left({X}^{+}_{{j}_p}(t)={X}^{+}_{{j}_A}(t)={X}^{+}_{{j}_B}(t)=0\right)
		=
		\mathbb{P} \left({X}^{+}_{{j}_p}(t)=0\right)
		\mathbb{P} \left({X}^{+}_{{j}_A}(t)={X}^{+}_{{j}_B}(t)=0\right).
	\end{equation}
  By Lemma~\ref{thm:[Si,Sj>=[Si][Sj]},
	\begin{equation}
		\label{eq:Prob(X^j,A)-conditional}
		\mathbb{P} \left({X}^{+}_{{j}_A}(t)={X}^{+}_{{j}_B}(t)=0\right)
		\ge 	\mathbb{P} \left({X}^{+}_{{j}_A}(t)=0 \right)
		\mathbb{P} \left({X}^{+}_{{j}_B}(t)=0 \right).
	\end{equation}
	Combining relations~\eqref{eq:Prob(X^j,A)}, 
	\eqref{eq:Prob(X^j,A)-indifference}, 
	and~\eqref{eq:Prob(X^j,A)-conditional}  gives
	\begin{equation}
		\label{eq:Prob(X^j,A)+indifference}
		\mathbb{P} \left(X_j(t)=0\right) 
		\ge
		\mathbb{P} \left({X}^{+}_{{j}_p}(t)=0\right)
		\mathbb{P} \left({X}^{+}_{{j}_A}(t)=0\right)
		\mathbb{P} \left({X}^{+}_{{j}_B}(t)=0\right).
	\end{equation}
	Substituting~\eqref{eq:Prob(X^j,A)=Prob(X^j,tilde-A)} in~\eqref{eq:Prob(X^j,A)+indifference} proves~\eqref{eq:funnel_inequality}. \Halmos
\endproof

In order to determine the conditions under which the funnel inequality becomes an equality, we  introduce the notion of a {\em funnel node}:

\begin{definition} [funnel node] \label{def:funnel} Let $\{A,B,j\}$ be a partition of~${\cal M}$.  Node~$j$ is called a ``funnel node of~$A$ and~$B$ in network~${\cal N}$'', if there is no node in~$A \cup B$ which is influential to~$j$ both in~${\cal N}^A$ and in~${\cal N}^B$.
\end{definition}

 Recall also the following definition:
 \begin{definition} [vertex  cut (vertex separator)]
  \label{def:vertex-cut} Let $\{A,B,j\}$ be a partition of~${\cal M}$.  Node~$j$ is called a  ``vertex  cut'' or ``vertex separator''
 between $A$ and $B$, if removing node~$j$ from the network makes the two sets~$A$ and~$B$ disconnected (see Figure~\ref{fig:funnel}A).
 \end{definition}

 Any node which is a vertex cut is also a funnel node:
\begin{lemma}
	  \label{lem:j-is-funnel}
	 Let $\{A,B,\{j\}\}$ be a partition of~${\cal M}$.
	 If node~$j$ is  a  vertex cut between~$A$ and~$B$,
 then $j$ is a  funnel node.
	\end{lemma}
\proof{Proof.}
	Let $m \in A$ be an influential node to~$j$. Then $m$ cannot be 
	an influential node to~$j$ in~${\cal N}^B$, since in~${\cal N}^B$
	we removed all edges from $A$ to~$j$, and so there is no sequence of edges (influential or not) from~$m$ to~$B$. \Halmos  
\endproof

The converse statement, however, is not true, i.e., 
there are networks in which $j$ is a funnel node, and~$A$ and~$B$ 
are directly connected. For example, this is the case if all edges
between~$A$ and~$B$ are non-influential to~$j$.  Moreover, even 
if nodes $a \in A$ and~$b \in B$ are connected by an influential edge
$\circled{a} \to \circled{b}$, $j$~may still be a funnel node, provided that there is no influential edge that emanate from node~$a$ in network~${\cal N}^A$ (e.g.\ Figure~\ref{fig:funnel}B).

\begin{theorem}[Funnel equality]
\label{thm:funnel_node_equality}
Consider the discrete Bass model~\eqref{eq:dbm}, and
let $\{A,B,j\}$ be a partition of~${\cal M}$. 
\begin{itemize}
	\item 

If $j$ is a funnel node of~$A$ and~$B$,  then~\eqref{eq:funnel_inequality} becomes 
 the {\bf funnel equality}  
\begin{equation}
	\label{eq:funnel_equality}
	\mathbb{P} \left(X^{A,B,p_j}_j(t)=0\right)  =
	\mathbb{P} \left(X^{A}_j(t)=0\right)
	\mathbb{P} \left(X^{B}_j(t)=0\right)
	\mathbb{P} \left(X^{p_j}_j(t)=0\right), \qquad t \ge0.
\end{equation}

\item 
If, however, $j$ is not a funnel node of~$A$ and~$B$, 
 then inequality~\eqref{eq:funnel_inequality} 
 is strict, i.e.,  
\begin{equation*}
	\label{eq:funnel_strong_inequality}
	\mathbb{P} \left(X^{A,B,p_j}_j(t)=0\right)  >
	\mathbb{P} \left(X^{A}_j(t)=0\right)
	\mathbb{P} \left(X^{B}_j(t)=0\right)
	\mathbb{P} \left(X^{p_j}_j(t)=0\right), \qquad t>0.
\end{equation*}
\end{itemize}
\end{theorem}
\proof{Proof.}
The inequality sign in the derivation of the funnel inequality~\eqref{eq:funnel_inequality} only comes from the use of Lemma~\ref{thm:[Si,Sj>=[Si][Sj]} in obtaining~\eqref{eq:Prob(X^j,A)-conditional}.  
By Lemma~\ref{thm:[Si,Sj>=[Si][Sj]}, inequality~\eqref{eq:Prob(X^j,A)-conditional} is
a strict inequality if and only if there exists a node~$m$ 
in network ${\cal N}^+$ which is influential to~${j}_A$ and to~${j}_B$.
Since no edges emanate from~${j}_A$, ${j}_B$, and~${j}_p$, then $m \in A \cup B$. 

Thus, the funnel inequality in strict if and only if there exists a node
$m \in A \cup B$ in network ${\cal N}^+$ which is influential to~${j}_A$ and to~${j}_B$.
This, however, is the case if and only if there exists a node $m \in A \cup B$
 which is influential to~${j}$ in network~${\cal N}^A$ and in~${\cal N}^B$,
 i.e.,  if~$j$ is not a funnel node of~$A$ and~$B$. \Halmos
\endproof
 The expressions~$\mathbb{P} \left(X^{A}_j(t)=0\right)$ and~$\mathbb{P} \left(X^{B}_j(t)=0\right)$ in Theorem~\ref{thm:funnel_node_equality} are usually unknown, as they do not take into account the effect of~$p_j$ on~$j$.
 Therefore, let us introduce two additional networks: On network ${\cal N}^{A,p_j}$, $j$ can adopt
 due to the combined influences of direct edges from~$A$ and due to~$p_j$, 
 and  on network ${\cal N}^{B,p_j}$, $j$ can adopt
 due to the combined influences of direct edges from~$B$ and due to~$p_j$.

\begin{definition} [Networks~${\cal N}^{A,p_j}$ and ${\cal N}^{B,p_j}$]
\label{def:more-networks}
Consider the discrete Bass model \eqref{eq:dbm} with network  
$
{\cal N} := {\cal N}({\cal M}, \{p_k\}_{k=1}^M, \{q_{k,j}\}_{k,j=1}^M)
$.
Let  $\{A,B,\{j\}\}$ be a partition of 
$ \cal{M}$. 
We define two additional networks with respect to this partition:
\begin{enumerate}
	\item ${\cal N}^{A,p_j}$  is obtained from~${\cal N}$ by removing  all direct links from~$B$ to~$j$, i.e., by setting~$q_{m,j} = 0$ for all~$m\in B$. Thus, we retain the direct edges  from~$A$ on~$j$ and the external influence~$p_j$. 
	\item ${\cal N}^{B,p_j}$ is defined similarly.
\end{enumerate}
\end{definition}

We can use the funnel equality to compute the combined influences 
from~$A$ and~$p_j$: 
\begin{lemma}
\label{lem:funnel_node_equality}
	Consider the discrete Bass model~\eqref{eq:dbm} on network ${\cal N}$.
Let $\{A,B,\{j\}\}$ be a partition of $ \cal{M}$.
Then
\begin{subequations}
	\label{eq:funnel_only2}
\begin{equation}
	\label{eq:funnel_only2A}
	\mathbb{P} \left(X^{A,p_j}_j(t)=0\right)  =
	\mathbb{P} \left(X^{A}_j(t)=0\right)
	\mathbb{P} \left(X^{p_j}_j(t)=0\right),  \qquad t \ge 0.
\end{equation}
\begin{equation}
	\label{eq:funnel_only2B}
	\mathbb{P} \left(X^{B,p_j}_j(t)=0\right)  =
	\mathbb{P} \left(X^{B}_j(t)=0\right)
	\mathbb{P} \left(X^{p_j}_j(t)=0\right),  \qquad t \ge 0.
\end{equation}
\end{subequations}
\end{lemma}
\proof{Proof.}
Let~$\widehat{\cal N}$ denote the network obtained from~${\cal N}^{A,p_j}$ by adding 
a fictitious isolated note, denoted by~$M+1$. Then
 $\{A \cup B, \{M+1\}, \{j\} \}$ is a partition of~$\{1, \dots, M+1\}$,
and  $j$~is a funnel node of~$A \cup B$ and
$\{M+1\}$ in~$\widehat{\cal N}$.
Let~$\widehat X_j$ denote the state of~$j$ in~$\widehat{\cal N}$. 
 By the funnel equality~\eqref{eq:funnel_equality},   
$$
\mathbb{P} \left(\widehat X^{A \cup B,\{M+1\} ,p_j}_j(t)=0\right)  =
\mathbb{P} \left(\widehat X^{A \cup B}_j(t)=0\right)
\mathbb{P} \left(\widehat X^{\{M+1\}}_j(t)=0\right)
\mathbb{P} \left(\widehat X^{p_j}_j(t)=0\right).
$$
By construction,
\begin{align*}
& \mathbb{P} \left(\widehat X^{A \cup B,\{M+1\} ,p_j}_j(t)=0\right)
 =\mathbb{P} \left(X^{A,p_j}_j(t)=0\right), 
 \qquad 
\mathbb{P} \left(\widehat X^{A \cup B}_j(t)=0\right)
=\mathbb{P} \left(X^{A}_j(t)=0\right), 
\\
&
\mathbb{P} \left(\widehat X^{\{M+1\}}_j(t)=0\right) \equiv 1, 
\qquad 
\mathbb{P} \left(\widehat X^{p_j}_j(t)=0\right)
=\mathbb{P} \left(X^{p_j}_j(t)=0\right),
\end{align*}
where~$X_j$ denote the state of~$j$ in network~${\cal N}$.
Therefore,  \eqref{eq:funnel_only2A}~follows. The proof for~\eqref{eq:funnel_only2B} is similar. \Halmos
\endproof

We can restate the funnel inequality and equality in terms of 
$\mathbb{P} \left(X^{A,p_j}_j(t)=0\right)$ and $\mathbb{P} \left(X^{B,p_j}_j(t)=0\right)$. 
This representation is useful when $\mathbb{P} \left(X^{A,p_j}_j(t)=0\right)$ and $\mathbb{P} \left(X^{B,p_j}_j(t)=0\right)$ correspond to known expressions
(see e.g., the proof of Theorem \ref{thm:alpha_q_circle}).

\begin{corollary}
\label{cor:funnel_node_inequality}
	Consider the discrete Bass model~\eqref{eq:dbm}.
Let $\{A,B,\{j\}\}$ be a partition of $ \cal{M}$.
Then 
\begin{equation}
	\label{eq:funnel_cor}
	\mathbb{P} \left(X^{A,B,p_j}_j(t)=0\right)  \ge
	\frac{\mathbb{P} \left(X^{A,p_j}_j(t)=0\right)
		\mathbb{P} \left(X^{B,p_j}_j(t)=0\right)}
	{\mathbb{P} \left(X^{p_j}_j(t)=0\right)},  \qquad t \ge 0.
\end{equation}
\begin{itemize}
	\item 
	
If $j$ is a funnel node of~$A$ and~$B$,  then inequality~\eqref{eq:funnel_cor} becomes an equality, i.e., 
	\begin{equation}
		\label{eq:funnel_equality-A-p}
		\mathbb{P} \left(X^{A,B,p_j}_j(t)=0\right)  =
		\frac{\mathbb{P} \left(X^{A,p_j}_j(t)=0\right)
			\mathbb{P} \left(X^{B,p_j}_j(t)=0\right)}
		{\mathbb{P} \left(X^{p_j}_j(t)=0\right)},  \qquad t \ge 0.
	\end{equation}
	
	\item 
	If, however, $j$ is not a funnel node of~$A$ and~$B$, then
inequality~\eqref{eq:funnel_cor} is strict, i.e., 
	\begin{equation}
		\label{eq:funnel_strict_inequality-A-p}
		\mathbb{P} \left(X^{A,B,p_j}_j(t)=0\right)  >
		\frac{\mathbb{P} \left(X^{A,p_j}_j(t)=0\right)
			\mathbb{P} \left(X^{B,p_j}_j(t)=0\right)}
		{\mathbb{P} \left(X^{p_j}_j(t)=0\right)},  \qquad t > 0.
	\end{equation}
\end{itemize}
\end{corollary}
\proof{Proof.}
This follows from Theorem~\ref{thm:funnel_node_inequality}, 
Theorem~\ref{thm:funnel_node_equality}, and Lemma~\ref{lem:funnel_node_equality}. \Halmos 
\endproof

\section{Circular networks.}
\label{subsec:homogeneous_circles}

We now present several applications of the funnel theorems. 
We begin with the discrete Bass model on a circle.
This problem was previously analyzed in~\cite{OR-10,Bass-boundary-18}.
In this section, we use the funnel theorem to derive a novel 
inequality  for diffusion on circles.

\subsection{Theory review.}

\begin{figure}[ht!]
	\begin{center}
		\scalebox{0.8}{\includegraphics{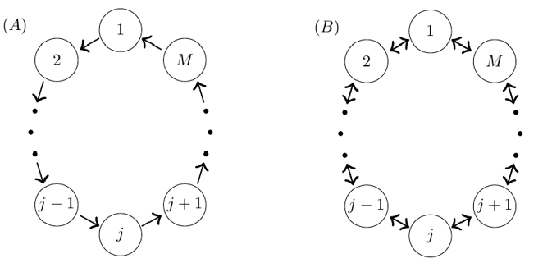}}
		\caption{(A)~One-sided circle. (B)~Two-sided circle.}
		\label{fig:circles}
	\end{center}
\end{figure}

We begin with a short theory review.
Let $f^{\rm 1-sided}_{\rm circle}(t;p,q,M)$ denote the expected fraction of adopters in a homogeneous one-sided circle with $M$ nodes, where each individual is only influenced by her left neighbor (see Figure~\ref{fig:circles}A), i.e., 
\begin{equation}
	\label{eq:p_j_q_j_one-sided-circle}
p_j\equiv p,\qquad q_{k,j}=\begin{cases}q, & \qquad {\rm if }\ (j-k)\Mod M=1,\\ 0, &\qquad {\rm if}\ (j-k) \Mod M\neq 1,\end{cases}\qquad j,k \in {\cal M}.
\end{equation}
Similarly, denote by $f^{\rm 2-sided}_{\rm circle}(t;p,q^R,q^L,M)$ the expected fraction of adopters in a homogeneous two-sided circle with
 $M$~nodes, where each node can be influenced by its left and right neighbors (see Figure~\ref{fig:circles}B), i.e.,
\begin{equation*}
	\label{eq:p_j_q_j_two-sided-circle}
p_j\equiv p,\qquad q_{k,j}=\begin{cases}
	q^L& \qquad {\rm if}\ (j-k)\Mod M=1,\\ q^R& \qquad  {\rm if}\ (k-j)\Mod M=1,\\  0&\qquad  {\rm if}\  |j-k|\Mod M\neq 1,
\end{cases}\qquad j,k \in {\cal M}.
\end{equation*}	

\begin{subequations}
	\label{eqs:compare_Fibich_Gibori}
	When
	\begin{equation}
		\label{eq:q_Lq_R}
		q=q^L+q^R,
	\end{equation}	
	the expected adoption fraction of adopters on one-sided and two-sided homogeneous circles is identical~\cite{OR-10}, i.e.,
	\begin{equation}
		\label{eq:compare_Fibich_Gibori}
		f_{\rm circle}^{{\rm 1-sided}}(t; p,q,M)\equiv f_{\rm circle}^{\rm 2-sided}(t;p,q^R,q^L,M), \qquad 0\leq t< \infty.
	\end{equation}
\end{subequations}
Therefore, we can drop the subscripts {\em 1-sided} and {\em 2-sided}.
For a finite circle, the expected fraction of adoption is given by the explicit expression~\cite[Lemma 4.1]{Bass-boundary-18}
\begin{subequations}
	\label{f_circle_computation}
	\begin{equation}
		f_{\rm circle}(t;p,q,M) = 1-\sum_{k=1}^{M-1} c_k
		\frac{\left(-q\right)^{k-1}}{p^{k-1}\left(k-1\right)!} e^{\left(-kp-q\right)t}  -  
		\frac{\left(-q\right)^{M-1}}{\prod_{j=1}^{M-1}\left(jp-q\right)} e^{-Mpt},
	\end{equation}
	where
	\begin{equation}
		\label{c_m-k}
		c_{M-k} =  1-\frac{q^k}{\prod_{j=1}^{k}\left(q-jp\right)}-\sum_{j=1}^{k-1}\frac{p^{j-k}\left(-q\right)^{-j+k}}{\left(k-j\right)!}c_{M-j}, \qquad k=1, \dots, M-1.
	\end{equation}
\end{subequations}
As $M\rightarrow \infty$ this expression simplifies into~\cite{OR-10}
\begin{equation}
	\label{eq:f_1D}
	\lim\limits_{M\to \infty}f_{\rm circle}(t;p,q,M)=f_{\rm 1D}(t;p,q):=1-e^{-(p+q)t+q\frac{1-e^{-pt}}{p}}.
\end{equation}

\subsection{ $[S_{\rm circle}](\cdot,q_1) \, [S_{\rm circle}](\cdot,q_2)<e^{-pt}[S_{\rm circle}](\cdot,q_1+q_2)$.}

\label{sec:S_circle^2<}
Let
 $$
 [S_{\rm circle}](t;p,q,M) := 1-f_{\rm circle}^{{\rm 1-sided}}(t;p,q,M)
$$
 denote the expected
	fraction of nonadopters in the Bass model~\eqref{eq:dbm} on the one-sided circle~\eqref{eq:p_j_q_j_one-sided-circle}.
We can use the funnel inequality to derive the following inequality, which is of interest by itself, and
 will also be used in the proofs of Theorems~\ref{thm:f_line^one-sided<f_line_two-sided} and~\ref{thm:f_1D_bound-periodic}:

\begin{theorem}
	\label{thm:alpha_q_circle}
	Let~$q_1,q_2>0$ and $3 \le M< \infty$. Then 
	\begin{equation}
		\label{eq:circle_convex}
		[S_{\rm circle}](t;p,q_1,M) \, [S_{\rm circle}](t;p,q_2,M)
		<e^{-pt} [S_{\rm circle}](t;p,q_1+q_2,M), \qquad t>0.
	\end{equation}
\end{theorem}
\proof{Proof.}
	Consider the Bass model on a two-sided anisotropic circle with~$M\geq 3$ nodes, 
	where the weights of the clockwise and counter-clockwise edges are~$q^R = q_1$ and~$q^L = q_2$, respectively. 
	Let~$\left\{j\right\}$, $A = \left\{1,\ldots, j-1\right\}$, and~$B: = \left\{j+1,\ldots,M\right\}$ be a partition of the nodes, where
	$1<j<M$.  The networks~${\cal N}^A$ and~${\cal N}^B$ are obtained from the original circle by removing the edges~$\circled{j} \leftarrow \circled{j+1}$ and~$\circled{j-1} \rightarrow \circled{j}$, respectively. Hence, any node~$m \in {\cal M} \setminus \{ j \}$ is influential to node~$j$ in both sub-networks~${\cal N}^A$ and~${\cal N}^B$. Consequently, 
	$j$~is {\em not}  a funnel node of~$A$ and~$B$. Therefore, the strict funnel inequality~\eqref{eq:funnel_strict_inequality-A-p} holds.
	
	The original network is a two-sided circle. Hence, by the equivalence of one-sided and two-sided circles, see~\eqref{eqs:compare_Fibich_Gibori},   
	\begin{subequations}
		\label{eqs:circle1-4}
		\begin{equation}
			\label{eq:circle1}
			\mathbb{P} \left(X_j^{A,B,p_j}(t)=  0\right)
			= [S_{\rm circle}](t;p,q_1+q_2,M).
		\end{equation}
		By~\eqref{eq:funnel_p_j},
		\begin{equation}
			\label{eq:circle2}
			\mathbb{P} \left(X_j^{p_j}(t)=  0\right) =
			e^{-pt}.
		\end{equation}
	
  	    The calculation of~$\mathbb{P} \left(X^{A, p_j}_j(t)=  0\right)$ is as follows.
	    In sub-network~${\cal N}^{A,p_j}$, we removed the edge~$\circled{j} \leftarrow \circled{j+1}$. As a result, all the clockwise edges $\{\circled{k} \leftarrow \circled{k+1}\}_{k \not=j}$ become noninfluential. Hence, by the indifference principle, 
	    we can compute $\mathbb{P} \left(X^{A, p_j}_j(t)=  0\right)$ on the
	    counter-clockwise one-sided circle with~$q^R=q_1$, i.e.,
		\begin{equation}
			\label{eq:circle3}
			\mathbb{P} \left(X_j^{A,p_j}(t)=  0\right) = [S_{\rm circle}](t;p,q_1,M).
		\end{equation}
		Using similar arguments, we can compute~$\mathbb{P} \left(X^{B,p_j}_j (t)=  0\right)$ on the equivalent clockwise one-sided circle with~$q^L=q_2$, yielding
		\begin{equation}
			\label{eq:circle4}
			\mathbb{P} \left(X_j^{B,p_j}(t)=  0\right) = [S_{\rm circle}](t;p,q_2,M).
		\end{equation}
	\end{subequations}
	Plugging relations~\eqref{eqs:circle1-4}
	into the funnel inequality~\eqref{eq:funnel_strict_inequality-A-p} gives the result. \Halmos
\endproof

\begin{remark}
 Since $e^{-pt} =   [S_{\rm circle}](t;p,0,M)$, 
inequality~\eqref{eq:circle_convex} can be rewritten as 
\begin{equation*}
	[S_{\rm circle}](t;p,q_1,M) \, [S_{\rm circle}](t;p,q_2,M)
	<  [S_{\rm circle}](t;p,0,M)  \,  [S_{\rm circle}](t;p,q_1+q_2,M).
\end{equation*}
\end{remark}

\begin{remark}
	 Substituting $q_1 = q_2 = \frac{q}{2}$ in~\eqref{eq:circle_convex} gives 
	$$[S_{\rm circle}]^2\left(t;p,\frac{q}{2},M\right) <  e^{-pt}  [S_{\rm circle}](t;p,q,M).
	$$
\end{remark}

The proof of Theorem~\ref{thm:alpha_q_circle} shows that 
inequality~\eqref{eq:circle_convex} is strict, because 
on the finite circle, $j$~is not a funnel node of 
$A = \left\{1,\ldots, j-1\right\}$ and~$B = \left\{j+1,\ldots,M\right\}$, 
and so the strict funnel inequality holds.
If we let $M \to \infty$, however,
removing node~$j$ from the network makes the two sets
$A = \left\{1,\ldots, j-1\right\}$ and~$B = \left\{j+1,\ldots,\infty \right\}$ 
disconnected.
Therefore,  on the infinite circle, $j$~is a funnel node of~$A$ and~$B$ 
(Lemma~\ref{lem:j-is-funnel}). Hence, on the infinite circle the funnel equality~\eqref{eq:funnel_equality-A-p} holds, and so inequality~\eqref{eq:circle_convex}
becomes an equality as~$M \to \infty$.  Indeed, we can prove this result directly:
\begin{lemma}
Let  $[S_{\rm 1D}]:=\lim_{M \to \infty} [S_{\rm circle}](t;p,q,M)$. Then
 	\begin{equation}
	\label{eq:S_1D_q1_q2}
	[S_{\rm 1D}](t;p,q_1) \, [S_{\rm 1D}](t;p,q_2)
	=e^{-pt} [S_{\rm 1D}](t;p,q_1+q_2), \qquad t \ge 0.
\end{equation}
\end{lemma}
\proof{Proof.}
 	Since $[S_{\rm 1D}]= e^{-(p+q)t+ q\frac{1-e^{-pt}}{p}}$, see~\eqref{eq:f_1D}, the result follows. \Halmos 
 \endproof

\section{Bounded lines.}

The discrete Bass model on a bounded line can be used to gain insight into the effects of boundaries on the diffusion. This problem was previously analyzed in~\cite{Bass-boundary-18}. In this section, we use the funnel theorems to 
obtain additional results.

\subsection{Bounded one-sided line.}
\label{sec:one-sided-line}

\begin{figure}[ht!]
\begin{center}
\scalebox{0.8}{\includegraphics{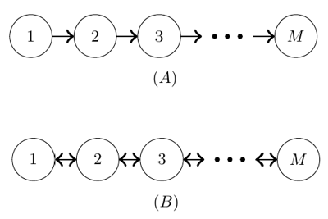}}
\caption{(A)~One-sided line. (B)~Two-sided line.}
\label{fig:lines}
\end{center}
\end{figure}

Consider the discrete Bass model~\eqref{eq:dbm} on the one-sided line~$[1,\ldots,M]$, where each node can only be influenced by its left neighbor (see Figure~\ref{fig:lines}A), i.e., 
	\begin{equation}
	\label{eq:p_j_q_j_onesided_line}
			p_j\equiv p,\qquad 
		q_{k,j}=
		\begin{cases}
			q, &  {\rm if }\ j-k=1,\\
			 0, &   \text{otherwise},
		\end{cases}
		\qquad  \qquad k,j\in {\cal M}.
	\end{equation}

 Let us denote the state of node~$j$ on the one-side line by $X_{j}^{\rm 1-sided}$. The adoption probability of node~$j$ on the one-sided line~\eqref{eq:p_j_q_j_onesided_line} is
 \begin{equation}
 	\label{eq:fj_one_sided_def}
 	 f_j^{\rm 1-sided}(t;p,q,M):=\mathbb{P} \left(X_{j}^{\rm 1-sided}(t)=1 \right).
 	 \end{equation}
   Unlike the Bass model on the circle, 
{\em there is no translation invariance on the line}, and so 
$f_j^{\rm 1-sided}$ may depend on~$j$.
Indeed, as $j$ increases, there are more nodes to its left that can adopt externally and then ``infect''~$j$ through a sequence of internal adoptions. 
Therefore, we have 
\begin{lemma}
	\label{lem:f_j^one-sided-monotone_in_j}
	$f_j^{\rm 1-sided}(t;p,q,M)$ is strictly monotonically-increasing  in~$j$.
\end{lemma}
\proof{Proof.}
By the indifference principle (Theorem~\ref{thm:Indifference}) with $\Omega = \{j\}$, 
\begin{equation}
	\label{eq:f_j_tilde_M}
  f_j^{\rm 1-sided}(t;p,q,M) \equiv f_j^{\rm 1-sided}(t;p,q,j), \qquad  M \ge j,
\end{equation}
since~$f_j$ does not change if we remove the non-influential edge
\circled{$j$} $\to$ \circled{$j+1$}.

Let ${\cal N}^A$ be the one-sided line~\eqref{eq:p_j_q_j_onesided_line} with $j \ge 2$ nodes, and let ${\cal N}^B$ denote the network obtained from~${\cal N}^A$ when we delete the influential edge \circled{$1$} $\to$ \circled{$2$}.
By the dominance principle (Lemma~\ref{lem:strong-dominance-principle-nodes}),
$$
  f_{j}^A(t) > f_{j}^B(t), \qquad 0<t<\infty.
$$  
 By~\eqref{eq:fj_one_sided_def} and~\eqref{eq:f_j_tilde_M},  
 $$
f_{j}^A(t)=f_j^{\rm 1-sided}(t;p,q,j)
=f_j^{\rm 1-sided}(t;p,q,M),
$$
 and 
$$
f_{j}^B(t)
=f_{j-1}^{\rm 1-sided}(t;p,q,M).
$$
 Therefore, the result follows. \Halmos
\endproof

Previously, we used the indifference principle to derive an explicit expression for~$f_j^{\rm 1-sided}$:
\begin{lemma}[\cite{Bass-boundary-18}]
	\label{lem:fj_one-sided_line}
		Let $f_{\rm circle}(t;p,q,j)$ denote the expected adoption level on a circle with $j$ nodes, see~\eqref{f_circle_computation}.  
%
	 Then
\begin{equation}
  \label{eq:f_j^onesided=f_circle(j)}
f_j^{\rm 1-sided}(t;p,q,M)= f_{\rm circle}(t;p,q,j). 
\end{equation}
\end{lemma}

From Lemmas~\ref{lem:f_j^one-sided-monotone_in_j} and~\ref{lem:fj_one-sided_line}, we have:
\begin{corollary}
	\label{cor:f_circle_increases_with_M}
	Let $t,p,q>0$. Then $f_{\rm circle}(t;p,q,M)$ is strictly monotonically-increasing  in~$M$.
\end{corollary}
Intuitively, fix  any node on a one-sided circle with $M$~nodes. There are $M-1$ nodes to its left that can adopt externally and than lead the node to adopt through a chain of internal adoptions. As $M$ increases, there are more nodes to its left. Hence, the adoption probability of the node increases.
Finally, for future reference, we note that since~$\lim_{M\rightarrow \infty}f_{\rm circle} = f_{\rm 1D}$, see~\eqref{eq:f_1D}, we have from Corollary~\ref{cor:f_circle_increases_with_M}  that 
	\begin{equation}
		\label{eq:f_circle<f1D}
		f_{\rm circle}(t;p,q,M)	<f_{\rm 1D}(t;p,q), \qquad 0<t, \quad M=1,2, \dots
	\end{equation}

\subsection{Bounded two-sided line.}
\label{sec:two-sided-line}
In~\cite{Bass-boundary-18}, we obtained an explicit but cumbersome expression for
the adoption probability of nodes on the two-sided line. In this section we use the funnel theorem to obtain a simpler expression.

Consider the discrete Bass model
on a two-sided homogeneous anisotropic line
with $M$~nodes, where each node can be influenced by its left and right neighbors
at the internal rates of~$q^L$  and~$q^R$, respectively  (see Figure~\ref{fig:lines}B).  
Thus,
	\begin{equation}
	\label{eq:p_j_q_j_twosided_line}
		p_j\equiv p,\qquad 
			q_{k,j}=
		\begin{cases}
			q^L,& ~~ {\rm if}\ j-k=+1,\\ 
			q^R,& ~~  {\rm if}\ j-k=-1,\\  
			0,& ~~		 \text{otherwise},
		\end{cases}
		\qquad 	\qquad k,j\in {\cal M}.
	\end{equation} 
Let us denote the adoption probability of node~$j$ in a two-sided line~\eqref{eq:p_j_q_j_twosided_line}  with
 $M$~nodes by
\begin{equation*}
	\label{eq:fj_two_sided_def}
	f_j^{\rm 2-sided}(t;p,q^L,q^R,M):=\mathbb{P} \left(X_{j}^{\rm 2-sided}(t)=1 \right),
\end{equation*}
and the corresponding nonadoption probability of node~$j$ by 
\begin{equation*}
	\label{eq:Sj_two_sided_def}
	[S_j^{\rm 2-sided}](t):=
	1-f_j^{\rm 2-sided}(t) = 
	\mathbb{P} \left(X_{j}^{\rm 2-sided}(t)=0 \right).
\end{equation*}   
  
In~\cite{fibich2022exact}, it was shown that $[S_j^{\rm 2-sided}]$
{\em can be written explicitly using $\left[S_j^{\rm 1-sided}\right]$
 for the one-sided left-going and right-going lines}, as follows.
Let~$[S_j^L](t)$ denoted the probability that $X_j(t)=0$ when 
we discard the influences of all the right neighbors
by setting $q^{R} \equiv0$ in~\eqref{eq:p_j_q_j_twosided_line}, so that the network becomes a left-going one-sided line. 
Similarly, we denote by $[S_j^R](t)$ the nonadoption probability
of node~$j$ on the right-going one-sided line, which is obtained  by setting $q^{L}\equiv0$ in~\eqref{eq:p_j_q_j_twosided_line}.

\begin{lemma} [\cite{fibich2022exact}]
	\label{lem:1D2SWithP_R-hom}
	Consider the discrete Bass model~\eqref{eq:dbm} on the two-sided line \eqref{eq:p_j_q_j_twosided_line}.
	Then 
	\begin{equation}
		\label{eq:EquationDiscrete1D2Sided_1}
		[S_j^{\rm 2-sided}](t)
		=
	\frac{[S_j^L](t)[S_j^R](t)}{e^{-p t}}, \qquad j\in\mathcal{M}, \qquad t \ge 0.
	\end{equation}
\end{lemma}

\proof{Proof.}
	 We provide a simpler proof for this lemma, which makes use of the funnel equality. Let $j$ be an interior node.  Let~$A := \left\{1,\ldots,j-1\right\}$ and~$B := \left\{j+1,\ldots,M\right\}$. Hence, $\{A,B,\{j\}\}$ is a partition of the nodes. 
	 Since the sets~$A$ and~$B$ become disconnected
	  if we remove node~$j$,  we have from Lemma~\ref{lem:j-is-funnel} that $j$~is a funnel node to~$A$ and~$B$. Therefore, by the funnel equality~\eqref{eq:funnel_equality-A-p}, 
	  \begin{equation*}
	  	\label{eq:funnel_thrm2-two-sided}
	  	\mathbb{P} \left(X^{A,B,p}_j\left(t\right)=0\right)=
	  	\frac{	\mathbb{P} \left(X^{A,p}_j\left(t\right)=0\right)\mathbb{P} \left(X^{B,p}_j\left(t\right)=0\right)}
	  	{\mathbb{P} \left(X^{p}_j\left(t\right)=0\right)}
	  	.
	  \end{equation*} 
	By the indifference principle, when we compute $[S_j]$ on a two-sided line, 
	we can delete all the noninfluential edges that point away from node~$j$. Therefore,
	$$
	\mathbb{P}\left(X^{A,p}_j\left(t\right)=0\right) =[S_j^L](t) , \qquad 
	\mathbb{P}\left(X^{B,p}_j\left(t\right)=0\right)=[S_j^R](t). 
	$$	
	In addition, 
	$
	\mathbb{P}\left(X^{p}_j\left(t\right)=0\right) = 
	e^{-pt}.
	$
	Hence, \eqref{eq:EquationDiscrete1D2Sided_1} follows.
	
	Let $j$ be the left boundary node. Then on the left-going one-sided line, 
	$j$~is not influenced by any node, and so $[S_j^L] = e^{-pt}$.
	In addition, by the indifference principle, 
	$[S_j^{\rm 2-sided}]=[S_j^R]$. 
	Therefore, \eqref{eq:EquationDiscrete1D2Sided_1} follows. The proof for the right boundary point is identical. \Halmos
\endproof

\subsubsection{Simpler expression for $f_j^{\rm 2-sided}$.}
An explicit expression for the adoption probability~$f_j^{\rm 2-sided}$ 
of nodes on a two-sided line with $M$~nodes
was previously obtained in~\cite{Bass-boundary-18} in the isotropic case~$q^L=q^R=\frac{q}{2}$. Thus,  $f_j^{\rm 2-sided} = 1-[S_j^{\rm 2-sided}]$,
where 
\begin{flalign*}
\label{eq:S_two_sided-old}
[S_j^{\rm 2-sided}](t;p,q^R=\frac{q}{2},q^L=\frac{q}{2},M) = 
\begin{cases}
[S_{\rm circle}](t;p,\frac{q}{2},M),\quad &\text{$j=1, M$,} \\
e^{-(p+q)t}\left(1+\frac{q}{2}A_j(t)\right), \quad &\text{$j=2, \dots, M-1$,}
\end{cases}
\end{flalign*}
where
\begin{equation*}
\label{eq:two_sided_line_s_tilde}
\begin{aligned}
A_j\left(t\right)& = \int_{0}^{t}e^{(p+q)\tau}\Big[[S_{\rm circle}](\tau;p,\frac{q}{2},j)[S_{\rm circle}](\tau;p,\frac{q}{2},M-j)\\
&\qquad \quad ~~+ [S_{\rm circle}](\tau;p,\frac{q}{2},j-1)[S_{\rm circle}](\tau;p,\frac{q}{2},M-j+1)\Big]d\tau.
\end{aligned}
\end{equation*}
A simpler expression, which is also valid in the anisotropic case
~$q^L\not=q^R$,  can be obtained using the funnel equality:
\begin{theorem}
	\label{thm:fj_twosided_new}
		Consider the discrete Bass model~\eqref{eq:dbm} on the two-sided line~\eqref{eq:p_j_q_j_twosided_line}.
	Then 
	\begin{equation}
		\label{eq:S_two_sided}
		[S_j^{\rm 2-sided}](t;p,q^R,q^L,M) =  e^{pt} \,[S_{\rm circle}](t;p,q^L,j) \, [S_{\rm circle}](t;p,q^R,M-j+1), \quad j \in {\cal M},
	\end{equation}
where $[S_{\rm circle}]:=1-f_{\text{\rm circle}}$ 
	and $f_{\rm circle}$ is given by~\eqref{f_circle_computation}.
\end{theorem}
\proof{Proof.}
By Lemma~\ref{lem:1D2SWithP_R-hom}, 
$$
	[S_j^{\rm 2-sided}](t)
	=
	\frac{[S_j^L](t)[S_j^R](t)}{e^{-p t}}.
$$
	On the two-sided line $[1, \dots, M]$,
$$
 [S_j^L](t) = [S_j^{\rm 1-sided}](t;p,q^L,M), \qquad 
 [S_j^R](t) = [S_{M-j+1}^{\rm 1-sided}](t;p,q^R,M),
  $$
  where the subscript~$k$ in~$[S_k^{\rm 1-sided}]$ refers to the node number,
  counted in the direction of the one-sided line. Since
  $
  [S_j^{\rm 1-sided}](t;p,q,M) = [S_{\rm circle}](t;p,q,j),
  $
  see~\eqref{eq:f_j^onesided=f_circle(j)}, the result follows. \Halmos
  \endproof

\subsubsection{$f_{\rm line}^{\rm 1-sided}(t;p,q,M)<f_{\rm line}^{\rm 2-sided}(t;p,q^L,q^R,M)$.} 
\label{sec:f_one-sided<f_two-sided}
As noted, when~$q = q^R+q^L$,
one-sided and two-sided diffusion on the circle are identical, i.e., 
 $f_{\rm circle}^{\rm 1-sided} \equiv f_{\rm circle}^{\rm 2-sided}$, see \eqref{eqs:compare_Fibich_Gibori}.
On finite lines, however, this is not the case. Indeed, in~\cite{Bass-boundary-18} we showed that one-sided diffusion is strictly
slower that isotropic two-sided diffusion 
i.e., 
$$
f_{\rm line}^{\rm 1-sided}(t;p,q,M)<
f_{\rm line}^{\rm 2-sided}\left(t;p,q^L= \frac{q}{2},q^R= \frac{q}{2},M\right), \qquad 0<t<\infty.
$$
The availability of the new explicit expression~\eqref{eq:S_two_sided}
for~$[S_j^{\rm 2-sided}]$ 
allows us to generalize this result to the anisotropic two-sided case
($q_R \neq q_L$), with a much simpler proof: 
\begin{theorem}
	\label{thm:f_line^one-sided<f_line_two-sided}
	Let  $q = q^L + q^R$. Then for any $p,q_L,q_R>0$ and $2 \le M <\infty$, 
	$$
	f_{\rm line}^{\rm 1-sided}(t;p,q,M)<f_{\rm line}^{\rm 2-sided}(t;p,q^L,q^R,M), \qquad 0<t<\infty.
	$$
\end{theorem}
\proof{Proof.}
	Let $t>0$.
	Since $f_{\rm line} = \sum_{j=1}^{M}f_{j}$,  
we need to prove that
\begin{equation*}
\frac1M \sum_{j=1}^{M}f^{\text{1-sided}}_{j}<	\frac1M \sum_{j=1}^{M} f^{\text{2-sided}}_{j}.   
\end{equation*}
The key to proving this inequality is to show that it holds for any {\em  pair of nodes $\{k,M+1-k\}$ which are symmetric about the midpoint},
i.e., that
\begin{equation}
	\label{eq:symmetric-pairs}
	f^{\rm 1-sided}_{k}+f^{\rm 1-sided}_{M+1-k}< f^{\rm 2-sided}_{k}+f^{\rm 2-sided}_{M+1-k}, \qquad k \in  {\cal M} .
\end{equation}
This inequality was originally proved in~\cite{Bass-boundary-18}.
That proof, however, was very long and technical. 
We now give a simpler proof of~\eqref{eq:symmetric-pairs}, which makes use of the
the new explicit expression~\eqref{eq:S_two_sided}, which was derived using the  funnel inequality.

Equation~\eqref{eq:symmetric-pairs} can be rewritten as   
$$
[S^{\rm 1-sided}_{k}]+[S^{\rm 1-sided}_{M+1-k}]> [S^{\rm 2-sided}_{k}]+[S^{\rm 2-sided}_{M+1-k}], \qquad k \in  {\cal M} .
$$
Therefore, by~\eqref{eq:f_j^onesided=f_circle(j)} and~\eqref{eq:S_two_sided}, it suffices to prove for $t>0$ that
\begin{align}
  \label{eq:apolo}
&[S_{\rm circle}]\left(t;p,q,k\right)+ 
[S_{\rm circle}]\left(t;p,q,M+1-k\right)
 >
 \\
 \nonumber
   &  \qquad\qquad\qquad \qquad\qquad
 e^{pt}[S_{\rm circle}]\left(t;p,q^L,k\right) [S_{\rm circle}]\left(t;p, q^R,M-k+1\right)
 \\
 \nonumber
 & \qquad\qquad\qquad\qquad\qquad +
  e^{pt}[S_{\rm circle}]\left(t;p, q^L,M-k+1\right)[S_{\rm circle}]\left(t;p, q^R,k\right) 
 .\end{align}
Let $s(q,k):=[S_{\rm circle}](t; p, q,k)$ and $\tilde{k}:=M-k+1$.
By Theorem~\ref{thm:alpha_q_circle}, 
$$
s(q,k) >  
e^{pt} s(q^R,k)s(q^L,k),  
\qquad s(q,\tilde{k}) > 
e^{pt}s(q^R,\tilde{k})s(q^L,\tilde{k}).
$$
Therefore, to prove~\eqref{eq:apolo}, it suffices to show that
\begin{align*}
&s(q^L,k) s(q^R,k)
+
s(q^L,\tilde{k})s(q^R,\tilde{k}) 
	>
	s(q^L,k) s(q^R,\tilde{k})
	 +
	s(q^L,\tilde{k})s(q^R,k), 
	\end{align*}
 i.e., that 
$$
\Big(s(q^L,k)-s(q^L,\tilde{k})\Big)
\Big(s(q^R,k)- s(q^R,\tilde{k})\Big)>0.
$$
This inequality
follows from the strict monotonicity of 
 $s(q,k):=[S_{\rm circle}]\left(t; p, q,k\right)$
 in~$k$, see Corrollary~\ref{cor:f_circle_increases_with_M}. Therefore, we proved~\eqref{eq:apolo}. \Halmos 
\endproof

The condition $q^R, q^L>0$ ensures that the two-sided line does not
trivially reduce to the one-sided line.

\section{$D$-dimensional Cartesian networks.}
  \label{sec:f_1D-lower-bound}
Consider an infinite $D$-dimensional homogeneous Cartesian network~$\mathbb{Z}^{\rm D}$, where nodes are labeled by their $D$-dimensional coordinate vector~${\bf j} = (j_1,\cdots,j_D) \in \mathbb{Z}^D$. For one-sided networks,  each node can be influenced by its $D$ nearest-neighbors at the rate of~$\frac{q}{D}$, and so the external and internal parameters are
\begin{subequations}
	\label{eq:Bass-model-D-dimensional-one-sided-eqs}
	\begin{equation}
		\label{eq:Bass-model-D-dimensional-one-sided}
		p_{\bf j} \equiv p, \qquad q_{{\bf k},{\bf j}} = \left\{ 
		\begin{array}{ll}
			\frac{q}{D}, & \quad {\bf j} - {\bf k} = {\bf \widehat{e}}_i ,\\
			0, & \quad  \mbox{otherwise},
		\end{array}
		\right. ,\qquad {\bf k}, {\bf j} \in \mathbb{Z}^D, \quad i \in \left\{1,\ldots, D\right\},
	\end{equation}
	where~${\bf \widehat{e}}_i \in \mathbb{Z}_{\rm D}$ is the unit vector in the $i$-th coordinate.  For two-sided Cartesian networks,
	 each node can be influenced by its $2D$ nearest-neighbors at the rate of~$\frac{q}{2D}$, and so the external and internal parameters are
	\begin{equation}
		\label{eq:Bass-model-D-dimensional}
		p_{\bf j} \equiv p, \qquad q_{{\bf k},{\bf j}} = \left\{ 
		\begin{array}{ll}
			\frac{q}{2D}, & \quad {\bf j} - {\bf k} = \pm {\bf \widehat{e}}_i ,\\
			0, & \quad  \mbox{otherwise},
		\end{array}
		\right. ,\qquad {\bf k}, {\bf j} \in \mathbb{Z}^D, \quad i \in \left\{1,\ldots, D\right\}. 
	\end{equation}
\end{subequations}
Note that for both~\eqref{eq:Bass-model-D-dimensional-one-sided} and~\eqref{eq:Bass-model-D-dimensional}, the weights of the edges are normalized so that
\begin{equation*}
	q_{\bf j} \equiv q, \qquad {\bf j}\in \mathbb{R}^D,
\end{equation*}
see~\eqref{eq:q_on_node}.
We denote the fraction of adopters  on one-sided and two-sided $D$-dimensional Cartesian networks by $f_D^{\rm 1-sided}$ 
and~$f_D$, respectively.

\subsection{$f_{\rm D}>f_{\rm 1D}$.}

\begin{figure}[ht!]
\begin{center}
\scalebox{0.7}{\includegraphics{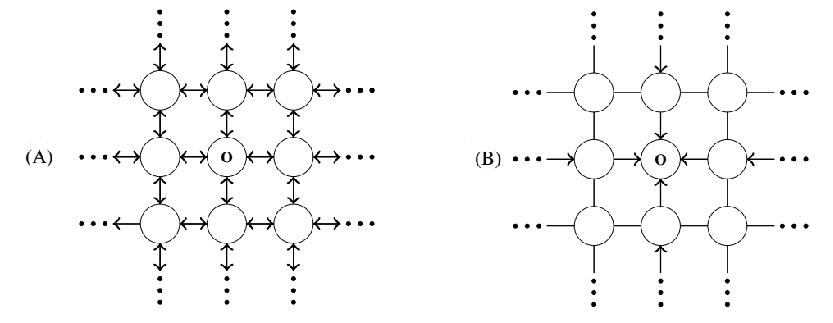}}
\caption{Illustration of networks used in the proof of Theorem~\ref{thm:f_1D_bound} is the two-sided $2D$ case. (A) ${\cal N}_{\rm 2D}$: Infinite 2D two-sided network. (B)~${\cal N}_{\rm 2D}^{\rm rays}$: After removing the edges from network~${\cal N}_D$, node~${\bf 0}:=\left(0,0\right)$ is at the intersection of~$2D=4$ one-sided rays.}
\label{fig:twodimension_funnel}
\end{center}
\end{figure}

In~\cite{OR-10}, it was observed numerically that 
$f_D$ is monotonically increasing in~$D$, i.e., 
$$
f_{\rm 1D}(t,p,q)<f_{\rm 2D}(t,p,q) < f_{\rm 3D}(t,p,q)< \cdots,   \qquad  t>0, 
$$ 
and similarly for one-sided diffusion. So far, however, 
this result was only proved for small times~\cite[Lemma~14]{fibich2022diffusion}.
We can use the funnel theorem to provide a partial proof, 
namely,  that
 $f_D>f_{1D}$ for all~$D\geq 2$:\,\footnote{Intuitively, this is because the diffusion evolves as a random creation of external seeds, which expands into clusters. The expansion rate of multi-dimensional clusters grows with the cluster size, whereas $1D$ clusters grow at a constant rate of~$q$. See~\cite{OR-10} for more details.} 
\begin{theorem} 
                \label{thm:f_1D_bound}
%
For any~$D\geq 2$ and $t,p,q>0$,
                \begin{equation*}
f_D(t;p,q) > f_{\rm 1D}(t;p,q), \qquad                                 
f_{\rm D}^{\rm 1-sided}(t;p,q)> f_{\rm 1D}(t;p,q),
                \end{equation*}
            where $f_{\rm 1D}$ is given by~\eqref{eq:f_1D}.
\end{theorem}
\proof{Proof.}
Let $D\geq 2$, and let~${\cal N}_D$ denote an infinite $D$-dimensional one-sided or two-sided Cartesian network. Denote 
the origin node by ${\bf 0}:=(0, \dots, 0) \in \mathbb{Z}^D$.
By  translation invariance,
\begin{equation}
                \label{eq:N1_cart}
                \mathbb{P}\left(X^{{\cal N}_D}_{{\bf k}={\bf 0}}(t) = 1\right) = f_{\rm D}\left(t;p,q\right).
\end{equation}
Let network~${\cal N}_D^{\rm rays}$ be obtained from~${\cal N}_D$ by removing all edges, except for those that lie on lines that go through the origin node~${\bf 0}$ and also point towards~${\bf 0}$ (see Figure~\ref{fig:twodimension_funnel}).
Hence, the origin node~${\bf 0}$ in network~${\cal N}_D^{\rm rays}$ is the intersection of~$D$ one-sided rays with edge weights~$\widetilde{q} = \frac{q}{D}$ in the one-sided network, and at the intersection of~$2D$ one-sided rays with edge weights~$\widetilde{q} = \frac{q}{2D}$ in the two-sided network. In Lemma~\ref{lem:intersection} below, we will prove that  
\begin{equation}
                \label{eq:N1_1D}
                \mathbb{P}\left(X^{{\cal N}_D^{\rm rays}}_{{\bf 0}}\left(t\right) = 1\right) =f_{\rm 1D}\left(t;p,q\right).
\end{equation}
Since some of the edges that were removed in~${\cal N}_D^{\rm rays}$ are influential to the origin node, then
by the strong dominance principle for nodes
(Lemma~\ref{lem:strong-dominance-principle-nodes}),
\begin{equation}
                \label{eq:N1>N2}
                \mathbb{P}\left(X^{N_D}_{{\bf 0}}\left(t\right) = 1\right)
                >\mathbb{P}\left(X^{N_D^{\rm rays}}_{{\bf 0}}\left(t\right) = 1\right).
\end{equation}
Combining~\eqref{eq:N1_cart}, \eqref{eq:N1_1D}, and \eqref{eq:N1>N2} gives the result. \Halmos
\endproof

\begin{figure}[ht!]
\begin{center}
\scalebox{0.8}{\includegraphics{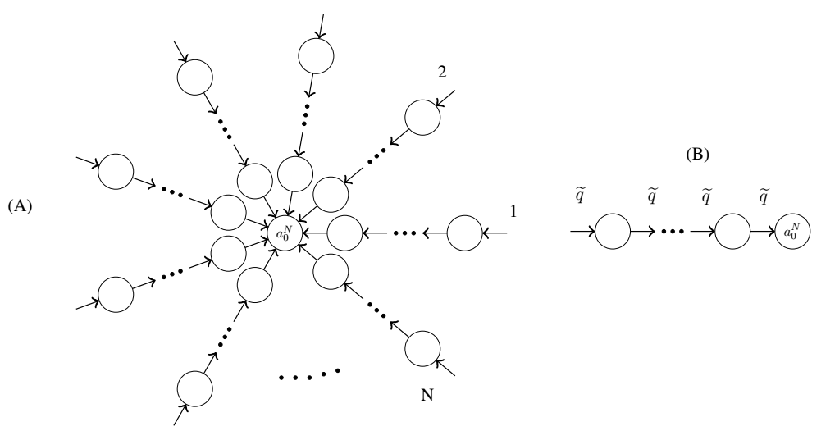}}
\caption{(A)~Node $a_0^N$ is at the intersection of~$N$ one-sided semi-infinite rays. (B)~A single semi-infinite ray. The weight of all edges is~$\widetilde{q}$.}
\label{fig:intersection}
\end{center}
\end{figure}

To finish the proof of Theorem~\ref{thm:f_1D_bound}, we prove
\begin{lemma}
                \label{lem:intersection}
                Let node~$a_0^{N}$ be the intersection of~$N$ identical one-sided semi-infinite rays, such that the weight of all edges is~$\widetilde{q}$ (see Figure~\ref{fig:intersection}). Then
                                \begin{equation}
                                \label{eq:intersection}
                                                \mathbb{P}\left(X_{a_0^{N}}(t) = 1\right) = f_{\rm 1D}(t;p,N\widetilde{q}).
                                \end{equation}
\end{lemma}

\proof{Proof.}
                We prove by induction on~$N$, the number of rays.
                By Lemma \ref{lem:fj_one-sided_line},
                               $f^{\rm 1-sided}_{j}(t,p,\tilde{q},M) = f_{\rm circle}(t,p,\tilde{q},j)$.
                               Therefore,
                \begin{equation}
                \label{eq:basecase}
                \mathbb{P}\left(X_{a_0^{N=1}}(t) = 1\right)
                =
                \lim_{M \to \infty} f^{\rm 1-sided}_{j=M}(t,p,\tilde{q},M)
                =             \lim_{M \to \infty}
                                f_{\rm circle}(t,p,\tilde{q},M)
                = f_{\rm 1D}(t;p,\widetilde{q}),
\end{equation}
see~\cite{Bass-boundary-18}, which is the induction base $N=1$.

                For the induction stage, we prove the equivalent result
                \begin{equation}
                                \label{eq:intersection-S}
                                \mathbb{P}\left(X_{a_0^{N}}(t) = 0\right) = [S_{\rm 1D}](t;p,N\widetilde{q}).
                \end{equation}
Thus, we assume that~\eqref{eq:intersection-S} holds for
                $N$~rays, and prove that it also holds for~$N+1$ rays.
                Indeed, let~$A$ denote the first~$N$ rays, and let~$B$ denote the $(N+1){\rm th}$~ray. Then node~$a_0^{N+1}$ is a funnel node of~$A$ and~$B$. Therefore, by funnel equality~\eqref{eq:funnel_equality-A-p},
                \begin{equation}
                                \label{eq:induction0}
                                \mathbb{P}\left(X_{a_0^{N+1}}(t) = 0\right )= \frac{\mathbb{P}\left(X^{A,p}_{a_0^{N+1}}(t)=0\right)\mathbb{P}\left(X^{B,p}_{{a_0^{N+1}}}(t)=0\right)}{\mathbb{P}\left(X^{p}_{{a_0^{N+1}}}(t)=0\right)}.
                \end{equation}
By construction, $\mathbb{P}\left(X^{A,p}_{a_0^{N+1}}(t)=0\right) = \mathbb{P}\left(X_{a_0^{N}}(t)=0\right)$.
                Hence, by the induction assumption~\eqref{eq:intersection-S},
                \begin{subequations}
                \label{eqs:induction123}
                \begin{equation}
                                \label{eq:induction1}
                                \mathbb{P}\left(X^{A,p}_{a_0^{N+1}}(t)=0\right) = [S_{\rm 1D}](t;p,N\widetilde{q}).
                \end{equation}
Similarly, by construction, $\mathbb{P}\left(X^{B,p}_{a_0^{N+1}}(t)=0\right) = \mathbb{P}\left(X_{a_0^{N=1}}(t)=0\right)$.
                Therefore, by \eqref{eq:basecase},         \begin{equation}
                                \label{eq:induction2}
                                \mathbb{P}\left(X^{B,p}_{a_0^{N+1}}(t)=0\right)
                                = [S_{\rm 1D}](t;p,\widetilde{q}).
                \end{equation}
                By~\eqref{eq:funnel_p_j},
                \begin{equation}
                                \label{eq:induction3}
                                \mathbb{P}\left(X^{p}_{{a_0^{N+1}}}(t)=0\right) = e^{-pt}.
                \end{equation}
                \end{subequations}
                Plugging relations~\eqref{eqs:induction123}
                into~\eqref{eq:induction0} and using~\eqref{eq:S_1D_q1_q2}
                gives
                \begin{equation*}
                                \begin{aligned}
                                                \mathbb{P}\left(X_{a_0^{N+1}}(t) = 0\right) = [S_{\rm 1D}](t;p,(N+1)\widetilde{q}),
                                \end{aligned}
                \end{equation*}
                as desired. \Halmos
\endproof

\subsection{Torodial networks.}

    We can also consider the one-sided and two-sided discrete Bass 
    models~\eqref{eq:Bass-model-D-dimensional-one-sided}
    and~\eqref{eq:Bass-model-D-dimensional}, respectively, on  
   the $D$-dimensional toroidal Cartesian network 
   $T_D:=[1, \dots, M_1]^D$, which is 
   periodic in each of the $D$~coordinates and has $M = M_1^D$ nodes.
  Let $f_{\rm T_D}^{\rm 1-sided}$ and~$f_{\rm T_D}$ denote the expected fraction of adopters on the one-sided and two-sided torus, respectively.
We can extend Theorem~\ref{thm:f_1D_bound} to $D$-dimensional toroidal networks,
as follows:
\begin{theorem}
                \label{thm:f_1D_bound-periodic}
                For any~$D\geq 2$, $M_1 \ge 2$ and $t,p,q>0$,
                \begin{equation*}
f_{\rm T_D}(t;p,q,M_1^D)>f_{\rm circle}(t;p,q, M_1),
\qquad 
  f_{\rm T_D}^{\rm 1-sided}(t;p,q,M_1^D)>f_{\rm circle}(t;p,q, M_1),
                \end{equation*}
                where~$f_{\rm circle}$ is given by~\eqref{f_circle_computation}.
\end{theorem}

\proof{Proof.}
                The only difference from the proof in the infinite-domain case is that in the last stage of the proof of Lemma~\ref{lem:intersection}, we use
                inequality~\eqref{eq:circle_convex} instead of equality~\eqref{eq:S_1D_q1_q2}. \Halmos
\endproof

\section{Discussion.}

  The analytic tools developed in this study have numerous applications,
  as already demonstrated in Sections~\ref{subsec:homogeneous_circles} --\ref{sec:f_1D-lower-bound} and in~\cite{bounds-23}.  
  Beyond the specific results of this study, it reveals the intricate relations between the inequality $[S_{i,j}] \ge [S_i][S_j]$, the funnel theorems, 
Chebyshev's integral inequalities, and the concepts of influential nodes
and funnel nodes.  While all of these results are new for the Bass model,
some results appeared in some form in the study of epidemiological models, 
as we describe next.

\subsection{Relation to epidemiological models.}
 \label{sec:epidemiological}

If we sets $p_j\equiv 0$ in the discrete Bass model~\eqref{eq:dbm}, 
we obtain the discrete Susceptible-Infected (SI) model on networks 
from epidemiology. Therefore, 
Theorems~\ref{thm:[Si,Sj>=[Si][Sj]}  and~\ref{thm:[Si,Sj>[Si][Sj]}
(for $[S_{i,j}] \ge [S_i][S_j]$),
and all the funnel theorems, hold for the discrete SI model. 

In~\cite{cator2014nodal}, Cator and Van Mieghem  proved that  $[S_{i,j}] \ge [S_i][S_j]$
in the SIS and SIR epidemiological models. In these models, infected individuals 
later recover, and recovered individuals can either become infected again (SIS) or are immune from getting infected again  (SIR).
That study, however, did not include the equivalent of our 
Theorem~\ref{thm:[Si,Sj>[Si][Sj]}, namely, the  conditions under which 
this inequality is strict and when it is actually an equality. Indeed, the role played by {\em influential nodes} is one of the methodological contribution of our study.

In~\cite{kiss2015exact}, Kiss et al.\ derived 
the funnel equality~\eqref{eq:funnel_equality} for the SIR model,
for nodes that are vertex cuts. Our funnel theorems are more general
  in two aspects. First, we show that an equality holds not only 
  when the node is a vertex cut, but also when the node is a funnel node
  which is not a vertex cut. Second, we show that when the node is not 
  a funnel node, a strict funnel inequality holds, 
  and we find the direction of the inequality.
  
  Finally, we note that the relation between the inequality 
  $[S_{i,j}] \ge [S_i][S_j]$ and the funnel theorems was not noted 
  in the above studies.

\subsection{Universal lower bounds.}

	In~\cite{OR-10}, Fibich and Gibori conjectured that the adoption level~$f(t)$
	on any infinite network that satisfies 
	$
	p_j \equiv p$ and $q_j \equiv q$, see~\eqref{eq:q_on_node}, 
	is bounded from below by that on the infinite circle,
	i.e., $f(t) \ge f_{\rm 1D}(t)$, see~\eqref{eq:f_1D}. In~\cite{bounds-23},  however, we proved that the optimal universal lower bound for~$f(t)$ for all finite and infinite networks that satisfy 
	$
	p_j \equiv p$ and $q_j \equiv q$ is given by  the adoption level on a two-node circle, 	i.e., 
	$f(t) \ge f_{\rm circle}(t;p,q,M=2)$.
	Since $f_{\rm circle}(t;p,q,M=2)< f_{\rm 1D}$, 
	see~\eqref{eq:f_circle<f1D}, this shows that $f_{\rm 1D}$
is {\em not} a  universal lower bound for~$f(t)$ for all networks.
Theorem~\ref{thm:f_1D_bound} shows, however, that $f_{\rm 1D}$ is a 
universal lower bound, for all 
one-sided and two-sided infinite multi-dimensional Cartesian networks.

 \begin{APPENDICES}
\section{Chebyshev's integral inequality.}
\subsection{Proof of Lemma~\ref{lem:chebyshev-1D-weights}.}
\label{app:chebyshev-1D-weights}
Since  $f$ and $g$ are both non-decreasing (or both non-increasing) functions, and~$w$ is positive,
$$
(f(x)-f(y)) (g(x)-g(y))w(x)w(y)  \ge 0, \qquad x,y \in [0, 1].
$$
Hence, since $\int_a^b w(x)\,dx = 1$,
\begin{equation*}
\begin{aligned}
	0 &\le \int_a^b \int_a^b (f(x)-f(y)) (g(x)-g(y))w(x)w(y) \, dx dy \\
	&= \int_a^b \int_a^b \left(f(x)g(x)+f(y)g(y)-f(x)g(y)-f(y)g(x)\right)w(x)w(y) \, dx dy \\
	&=2  \int_a^b f(x) g(x)w(x) \, dx -  2  \int_a^b f(x)w(x)\, dx \int_a^b g(x)w(x)\, dx.
\end{aligned}
\end{equation*}

An equality holds if and only if   $(f(x)-f(y)) (g(x)-g(y))  \equiv 0$ almost everywhere in~$[a, b]^2$, which is the case  if and only if   either~$f(x)$ or~$g(x)$ are constants. \Halmos

\subsection{Proof of Lemma~\ref{lem:chebyshev-multi-D}.}
\label{app:chebyshev-multi-D} 
We prove by induction on~$D$. The case  $D=1$
is Lemma~\ref{lem:chebyshev-1D-weights} with~$[a,b]=[0,1]$ and~$w(x)\equiv 1$. 
For the induction stage, we assume that when~$f\left(x_1,\ldots,x_D\right)$ and $g\left(x_1,\ldots,x_D\right)$ are both non-decreasing (or both non-increasing) with respect to each~$x_i$ for $i=1,\ldots, D$, then
\begin{equation*}
	\begin{aligned}
		\int_{\left[0,1\right]^D}f&\left(x_1,\ldots,x_D\right)g\left(x_1,\ldots,x_D\right)dx_1\cdots dx_D \geq \\
		&\int_{\left[0,1\right]^D}f\left(x_1,\ldots,x_D\right)dx_1\cdots dx_D\int_{\left[0,1\right]^D}g\left(x_1,\ldots,x_D\right)dx_1\cdots dx_D,
	\end{aligned}
\end{equation*}
and  an equality holds if and only if for  any $i=1, \dots D$, 
either~$f$ or~$g$ are independent of~$x_i$.
We note that if $f\left(x_1,\ldots,x_{D+1}\right)$ is monotone  with respect to each~$x_i$ for $i=1,\ldots, D+1$, then $\int_{[0,1]}f\left(x_1,\ldots,x_{D+1}\right) dx_{D+1}$ is monotone with respect to each~$x_i$ for $i=1,\ldots, D$. 
Therefore, when $f\left(x_1,\ldots,x_{D+1}\right)$ and $g\left(x_1,\ldots,x_{D+1}\right)$ are both non-decreasing (or both non-increasing) with respect to each~$x_i$ for $i=1,\ldots, D+1$, then
\begin{equation*}
	\begin{aligned}
		&	\int_{\left[0,1\right]^{D+1}}f \left(x_1,\ldots,x_{D+1}\right)g\left(x_1,\ldots,x_{D+1}\right)dx_1\cdots dx_{D+1} \\
		& \qquad =\int_{\left[0,1\right]^{D}} \left[\int_{\left[0,1\right]}f\left(x_1,\ldots,x_{D+1}\right)g\left(x_1,\ldots,x_{D+1}\right)dx_{D+1}\right] dx_1\cdots dx_{D}\\
		&\qquad \geq \int_{\left[0,1\right]^{D}} \left[\int_{\left[0,1\right]}f\left(x_1,\ldots,x_{D+1}\right)dx_{D+1}\right] \left[\int_{\left[0,1\right]} g\left(x_1,\ldots,x_{D+1}\right)dx_{D+1}\right]dx_1\cdots dx_{D}\\
		& \qquad \geq \int_{\left[0,1\right]^{D}}\int_{\left[0,1\right]}f\left(x_1,\ldots,x_{D+1}\right)dx_{D+1}dx_1\cdots dx_D \int_{\left[0,1\right]^{D}}\int_{\left[0,1\right]}g\left(x_1,\ldots,x_{D+1}\right)dx_{D+1}dx_1\cdots dx_D \\
		& \qquad = \int_{\left[0,1\right]^{D+1}}f\left(x_1,\ldots,x_{D+1}\right)dx_1\cdots dx_{D+1}\int_{\left[0,1\right]^{D+1}}g\left(x_1,\ldots,x_{D+1}\right)dx_1\cdots dx_{D+1},
	\end{aligned}
\end{equation*}
where the first inequality follows from the monotonicity in~$x_{D+1}$
and the induction base, and second inequality follows from the
the monotonicity in~$\{x_{1}, \dots, x_{D}\}$ and the induction assumption.
Therefore, we proved the induction stage.

Finally, an equality holds if both of the inequalities above are equalities.  By Lemma~\ref{lem:chebyshev-1D-weights}, the first inequality is an equality if and only if   either $f$ or $g$ are independent of $x_{D+1}$. 
By the induction assumption, the second inequality is an equality if and only if  
for any $i=1, \dots, D$, either $\int_{\left[0,1\right]}f\left(x_1,\ldots,x_{D+1}\right)dx_{D+1}$ or $\int_{\left[0,1\right]}g\left(x_1,\ldots,x_{D+1}\right)dx_{D+1}$ are independent of~$x_i$,
which is the case if and only if   for any $i=1, \dots D$, either~$f$ or~$g$ are independent of $x_i$. \Halmos

	\section{Proof of Lemma~\ref{lem:fi,j-fifj=Si,j-SiSj}.}
\label{app:fi,j-fifj=Si,j-SiSj}
Since
\begin{align*}
	f_{i,j} & = \mathbb{P}(X_i(t) =1)-\mathbb{P}(X_i(t) =1,  X_j(t) = 0)
	\\ & = 
	\mathbb{P}(X_i(t) =1)-\mathbb{P}(X_j(t) =0) +\mathbb{P}(X_i(t) =0,  X_j(t) = 0) = f_i-[S_j]+[S_{i,j}],
\end{align*}
and 
$f_i  = 1-[S_i]$, then
\begin{align*}
	f_{i,j}-[S_{i,j}]   = f_i-[S_j] = 1-[S_i]-[S_j]. 
\end{align*}
In addition,
$$
f_i f_j-[S_i]\, [S_j] = (1-[S_i])(1-[S_j]) - [S_i]\, [S_j]  = 1-[S_i]-[S_j].
$$
Therefore, the result follows. \Halmos

	\section{Proof of~\eqref{eq:Prob(X^j,A)}.}
	\label{app:ProofX^jA}	
	   Let us fix $t>0$. Let us consider the Bass model on 
	   networks~${\cal N}$ and~${\cal N}^+$ with discrete times
	   $$
	  t_n := n \Delta t, \quad  n=0,1,\dots
	  $$ 
	  where~$\Delta t  = \frac{t}{N}$ and $N \gg 1$. 
	  Note that as $N \to \infty$, $\Delta t \to 0$ and $t_N \equiv t$.
	  It is thus sufficient to prove that 
	\begin{equation}
		 \label{eq:Prob(X^j,A)-discrete}
		\lim_{N \to \infty} \mathbb{P} \left(X_j(t_N)=0\right) =\lim_{N \to \infty}
		\mathbb{P} \left({X}^{+}_{{j}_p}(t_N)
		={X}^{+}_{{j}_A}(t_N)
		={X}^{+}_{{j}_B}(t_N)
		=0\right).
	\end{equation}

	To prove~\eqref{eq:Prob(X^j,A)-discrete}, we introduce the time-discrete realizations, see \eqref{eq:tilde-Xj^A-def},
$$
		\widetilde{X}_{k}(t_{n})
	:=X_{k}(t_{n};\{\vomega^{n}\}_{n=1}^\infty,\Delta t), \quad	\widetilde{X}_{k}^{+}(t_{n})
		:=X_{k}^{+}(t_{n};\{\vomega^{n}\}_{n=1}^\infty,\Delta t), \qquad  k \in {\cal M}, \quad  n = 0,1, \dots
$$
  	We also define the sub-realization $\{\vomega_{-j}^n\}_{n=1}^{\infty}$, where 
	$\vomega_{-j}^n := \{\omega_k^n\}_{k \in A \cup B}$. 
	Since there are no edges emanating from $j$, ${j}_A$, 
	${j}_B$, and $ {j}_p$,
	this sub-realization completely determines $\{\widetilde{X}_{k}(t_{n})\}$ and  $\{\widetilde{X}_{k}^+(t_{n})\}$ for all~$k\not=j$ and~$n \in \mathbb{N}$.
	Moreover, if we use the same $\{\vomega_{-j}^n\}_{n=1}^{\infty}$
	and~$\Delta t$ for both networks, then
	\begin{equation}
		\label{eq:tildeX_k(t_n)=X_k(t_N}
		\widetilde{X}_{k}(t_{n})
		 \equiv
		\widetilde{X}_{k}^+(t_{n}), \qquad 
		k \in A\cup B, \quad n = 0,1, \dots 
	\end{equation}

	To compute the left-hand-side of~\eqref{eq:Prob(X^j,A)-discrete},
	we first note that  
	$$
	\widetilde{X}_j(t_N)=0 \quad \iff \quad 
	 {\omega}_j^n \ge  \left(p_j+\sum\limits_{k \in A\cup B}  q_{k,j} \widetilde{X}_{k}(t_{n-1})\right) \Delta t, \quad n=1, \dots, N.
	$$
	Therefore, 
	$$
	\mathbb{P} \left(X_j(t)=0  \,\big |\, \{\vomega_{-j}^n\}_{n=1}^{N} \right)
	=  \prod_{n=1}^{N} F_n,   
	\qquad 
      F_n:=	
	1-\left(p_j+\sum_{k \in A \cup B}  q_{k,j} 
	\widetilde{X}_{k}(t_{n-1})
	\right) \Delta t.
	$$
	Hence, 
	\begin{align}
		\label{eq:P(X^{A,B,p}_j(t)=0)}
	&	\mathbb{P} \left(X_j(t)=0\right)
		= 
		\int_{[0, 1]^{(M-1)\times N}} \left( \prod_{n=1}^{N}  F_n\left(\{\vomega_{-j}^n\}_{n=1}^{N}\right) \right)
		d\vomega_{-j}^{1}\cdots d\vomega_{-j}^{N},
	\end{align}

	Similarly, to compute the right-hand-side
	 of~\eqref{eq:Prob(X^j,A)-discrete} , we note that
	 $$
	 \widetilde{X}^{+}_{{j}_p}(t_N)
	 =\widetilde{X}^{+}_{{j}_A}(t_N)
	 =\widetilde{X}^{+}_{{j}_B}(t_N)
	 =0
	 $$
	  if and only if  
	 $$
	 {\omega}_{{j}_p}^n \ge p_j \Delta t, \quad 
	{\omega}_{{j}_A}^n \ge \left(\sum\limits_{k \in A}  q_{k,j} \widetilde{X}^{+}_k(t_{n-1}) \right) \Delta t,
	\quad 
	{\omega}_{{j}_B}^n \ge \left(\sum\limits_{k \in B}  q_{k,j} \widetilde{X}^{+}_k(t_{n-1}) \right) \Delta t, \quad  n = 1, \dots, N.
	$$
	Therefore, 
	\begin{equation}
		\label{eq:P(tildeX^{A,p}_j(t)=0)}
		\mathbb{P} \left({X}^{+}_{{j}_p}(t) 
		={X}^{+}_{{j}_A}(t)
		={X}^{+}_{{j}_B}(t)
		=0\right)
		= 
	\int_{[0, 1]^{(M-1)\times N}}
	\left( \prod_{n=1}^{N}  F_n^+\left(\{\vomega_{-j}^n\}_{n=1}^{N}\right) \right)
	  d\vomega_{-j}^{1}\cdots d\vomega_{-j}^{N},
	\end{equation}
where
$$
F_n^+:= (1-p_j\Delta t)\left(1-\left(\sum_{k \in A}  q_{k,j} \widetilde{X}^{+}_k(t_{n-1}) \right) \Delta t\right) \left(1-\left(\sum_{k \in B}  q_{k,j} \widetilde{X}^{+}_k(t_{n-1}) \right) \Delta t\right) .
$$

   To finish the proof of~\eqref{eq:Prob(X^j,A)-discrete}, 
  we now show that the integrand $\prod_{n=1}^{N}	
  F_n$ of~\eqref{eq:P(X^{A,B,p}_j(t)=0)} approaches, uniformly in~$\{\vomega_{-j}^n\}_{n=1}^{N}$, the integrand $\prod_{n=1}^{N}	
  F_n^+$ of~\eqref{eq:P(tildeX^{A,p}_j(t)=0)}. Indeed,  
  by~\eqref{eq:tildeX_k(t_n)=X_k(t_N}, 
	\begin{eqnarray*}
		F_n^+&=&
		\left(1 - p_j\Delta t\right) \left(1-\left(\sum_{k \in A}  q_{k,j} \widetilde{X}^{+}_k(t_{n-1}) \right) \Delta t\right)
		\left(1-\left(\sum_{k \in B}  q_{k,j} \widetilde{X}^{+}_k(t_{n-1}) \right) \Delta t\right)
		\\
		&=&
\left(1 - p_j\Delta t\right) \left(1-\left(\sum_{k \in A}  q_{k,j} \widetilde{X}_k(t_{n-1}) \right) \Delta t\right)
\left(1-\left(\sum_{k \in B}  q_{k,j} \widetilde{X}_k(t_{n-1}) \right) \Delta t\right)
\\
			&=& 
		 \left(1-\left(p_j+\sum_{k \in A \cup B}  q_{k,j} \widetilde{X}_k(t_{n-1}) \right) \Delta t\right)+ A_n (\Delta t)^2 
		\\
		&=& 
		\left(1-\left(p_j+\sum_{k \in A \cup B}  q_{k,j} \widetilde{X}_k(t_{n-1}) \right) \Delta t\right) 
		\left(1+ \frac{A_n (\Delta t)^2}{1-\left(p_j+\sum_{k \in A \cup B}  q_{k,j} X_k(t_{n-1}) \right) \Delta t} \right)
		\\
	&	=&F_n	\left(1+ \frac{A_n (\Delta t)^2}{1-\left(p_j+\sum_{k \in A \cup B}  q_{k,j} X_k(t_{n-1}) \right) \Delta t} \right),
	\end{eqnarray*}
	where
	$$
	A_n = p_j \left(\sum_{k \in A \cup B}  q_{k,j} \widetilde{X}_k(t_{n-1}) \right)  +
	\left(\sum_{k \in A}  q_{k,j} \widetilde{X}_k(t_{n-1}) \right)\left(\sum_{k \in B}  q_{k,j} \widetilde{X}_k(t_{n-1}) \right)
	(1- p_j \Delta t).
	$$
	Hence,
	\begin{eqnarray}
		\label{eq:three-products}
		&&
		\prod_{n=1}^{N}	
		F_n^+
		= \prod_{n=1}^{N} 
		                     F_n
		\prod_{n=1}^{N} \left(1+ \frac{A_n(\Delta t)^2}{1-\left(p_j+\sum_{k \in A \cup B}  q_{k,j} X_k(t_{n-1}) \right) \Delta t} \right).
		\end{eqnarray}
By~\eqref{eq:P(X^{A,B,p}_j(t)=0)}--\eqref{eq:three-products}, to finish the proof of~\eqref{eq:Prob(X^j,A)-discrete}, we need to show that
	\begin{equation}
		\label{eq:prod->1-uniform}
		\lim_{\Delta t \to 0} \prod_{n=1}^{N} \left(1+ \frac{A_n (\Delta t)^2}{1-\left(p_j+\sum_{k \in A \cup B}  q_{k,j} X_k(t_{n-1}) \right) \Delta t} \right)  = 1,
	\end{equation}
	uniformly in $\{\vomega_{-j}^n\}_{n=1}^{N}$.
	The three sums that appear in $A_n$ are uniformly bounded:
	$$
	0 \le 
	\sum_{k \in A}  q_{k,j} X_k(t_{n-1}), \, \sum_{k \in B}  q_{k,j} X_k(t_{n-1})   \le 
	\sum_{k \in A \cup B}  q_{k,j} X_k(t_{n-1}) 
	\le \sum_{k \not= j}  q_{k,j}  = q_j. 
	$$
	Hence, for $0<\Delta t \ll 1$,
	$$
	0 \le A_n \le q_j\left(p_j +q_j\right) ,
	\quad
	   1 \le  \frac{1}{1-\left(p_j+\sum\limits_{k \in A \cup B}  q_{k,j} X_k(t_{n-1}) \right) \Delta t}  \le \frac1{1-(p_j+q_j) \Delta t} \le 2,
	$$
	and so each term in the product~\eqref{eq:prod->1-uniform} is bounded by 
	$$
	1 \le
	1+ \frac{A_n (\Delta t)^2}{1-\left(p_j+\sum\limits_{k \in A \cup B}  q_{k,j} X_k(t_{n-1}) \right) \Delta t} 
	\le 1+2q_j\left(p_j +q_j\right) (\Delta t)^2.
	$$
	Therefore, 
	$$
	1 \le 	\prod_{n=1}^{N} \left(1+ \frac{A_n (\Delta t)^2}{1-\left(p_j+\sum_{k \in A \cup B}  q_{k,j} X_k(t_{n-1}) \right) \Delta t} \right) \le \Big(1+2q_j\left(p_j +q_j\right) (\Delta t)^2 \Big)^{t/\Delta t}  . 
	$$
	As $\Delta t \to 0$, the right-hand-side approaches~1,  uniformly in $\{\vomega_{-j}^n\}_{n=1}^{N}$. Hence, we proved~\eqref{eq:prod->1-uniform}.
		Relation~\eqref{eq:Prob(X^j,A)} follows from~\eqref{eq:P(X^{A,B,p}_j(t)=0)},  \eqref{eq:P(tildeX^{A,p}_j(t)=0)}, \eqref{eq:three-products},
	and~\eqref{eq:prod->1-uniform}.
	Finally, since $j$ is an isolated node in network~${\cal N}^{p_j}$, we have~\eqref{eq:funnel_p_j}. \Halmos
 \end{APPENDICES}






\end{document}